\documentclass[10pt]{article}

\usepackage{amsmath,amsthm,amssymb}
\usepackage{caption,graphicx,subfig}
\usepackage[text={6.5in,9.25in},centering]{geometry}
\usepackage[sort,numbers]{natbib}
\usepackage{url}

\usepackage{xcolor}

\setcaptionmargin{0.25in}

\makeatletter\@addtoreset{equation}{section}\makeatother

\newtheoremstyle{named}{}{}{\itshape}{}{\bfseries}{.}{.5em}{\thmnote{#3's }#1}
\theoremstyle{named}

\theoremstyle{definition}

\newcommand{\Xv}{\mathbf{X}}

\newcommand{\fv}{\mathbf{f}}

\newcommand{\xv}{\mathbf{x}}

\newcommand{\Thetav}{\boldsymbol{\Theta}}

\newcommand{\Xiv}{\boldsymbol{\Xi}}
\newcommand{\Piv}{\boldsymbol{\Pi}}
\newcommand{\thetav}{\boldsymbol{\theta}}
\newcommand{\xiv}{\boldsymbol{\xi}}
\newcommand\Rb{\mathbb{R}}

\begin{document}

\title{Poincar\'e Maps for Multiscale Physics Discovery and Nonlinear Floquet Theory}

\author{
Jason J. Bramburger\thanks{Division of Applied Mathematics, Brown University, Providence, RI, 02906}\ \thanks{Department of Mathematics and Statistics, University of Victoria, Victoria, BC, V8P 5C2} \and 
J. Nathan Kutz\thanks{Department of Applied Mathematics, University of Washington, Seattle, WA, 98195-3925}
}

\date{}
\maketitle

\begin{abstract}
Poincar\'e maps are an integral aspect to our understanding and analysis of nonlinear dynamical systems. Despite this fact, the construction of these maps remains elusive and is primarily left to simple motivating examples. In this manuscript we propose a method of data-driven discovery of Poincar\'e maps based upon sparse regression techniques, specifically the sparse identification of nonlinear dynamics (SINDy) algorithm. This work can be used to determine the dynamics on and near invariant manifolds of a given dynamical system, as well as provide long-time forecasting of the coarse-grained dynamics of multiscale systems.  Moreover, the method provides a mathematical formalism for determining nonlinear Floquet theory for the stability of nonlinear periodic orbits. The methods are applied to a range of examples including both ordinary and partial differential equations that exhibit periodic, quasi-periodic, and chaotic behavior.          
\end{abstract}


\section{Introduction}\label{sec:Intro} 

At the turn of the 20th century, Henri Poincar\'e laid the foundations of what would become modern dynamical systems~\cite{Poincare}. It was particularly his consideration of the three-body problem that led to our first understanding of deterministic chaos. He also proposed the first recurrence map, or Poincar\'e map, characterizing the intersection of a periodic orbit in the state space of a continuous dynamical systems with a lower-dimensional, and transverse to the flow, subspace called the Poincar\'e section. The Poincar\'e map produces a discrete dynamical system with a state space one dimension lower than the original continuous dynamical system. The dynamics on the Poincar\'e map preserves many of the periodic and quasi-periodic orbits of the original system, and due to its dimensionality reduced form, it is often simpler to analyze than the original system.  This was especially true for celestial mechanics which was ideally suited for analysis via Poincar\'e maps~\cite{Poincare}.    

In practice Poincar\'e maps can be difficult to define explicitly as there are no general methods for their construction, but the advent of data-driven modeling techniques provide a suite of regression techniques, broadly falling under the aegis of machine learning, that can discover and/or fit dynamical models to time-series measurements of a dynamical system~\cite{BK2019}.  Principled methods for modeling multiscale dynamical systems must overcome the unique challenge of disambiguating between distinct spatial and/or temporal scales and Poincar\'e maps are a foundational aspect of dynamical systems that can be exploited for achieving this task and discovering multiscale physics.  Indeed, recent sampling strategies have already been shown to be effective for extracting multiscale temporal dynamics and the underlying micro-to-macro scale couplings~\cite{Champion18}. Poincar\'e maps provide a timescale separation by producing snapshots of the dynamics on the sampling scale of the map period. The recently developed {\em sparse identification of nonlinear dynamics} (SINDy)~\cite{SINDy} algorithm can then be used to discover governing equations of the map dynamics themselves, thus giving a macro-scale timescale description of the nonlinear dynamics. Therefore, in this manuscript we demonstrate the integration of SINDy and Poincar\'e maps on a diverse set of example problems from ordinary and partial differential equations and show that the decomposition provides a robust method for the discovery of multiscale nonlinear dynamics.
 
Our objective is to leverage the dimensionality reduction afforded by the Poincar\'e map in order to discover the nonlinear dynamics of the return map, or more precisely, the nonlinear dynamics of the slow scale dynamics. Linear models of the Poincar\'e dynamics can be constructed using Floquet theory~\cite{Fenichel,Hale,Kuehn} which is typically used to characterize the stability of periodic solutions. Floquet theory can alternatively be interpreted as giving linear, exponential solutions for the dynamics, much like the dynamic mode decomposition~\cite{Kutz2018}.  Instead of pursuing such linear methods, we modify the SINDy algorithm to discover the nonlinear dynamics that characterize the return map, thus providing a generalization of Floquet theory, or  {\em nonlinear Floquet theory}, and a more accurate representation of the underlying nonlinear dynamics. Indeed, the method can be used as a general architecture to evaluate the nonlinear stability of periodic orbits and invariant manifolds which go far beyond the simple perturbative frameworks of traditional dynamical systems analysis. Here we are able to demonstrate the integration of SINDy, Poincar\'e maps and Floquet theory on a diverse range of examples. The method highlights how the Poincar\'e map can be used to disambiguate and discover multiscale temporal dynamics, specifically the slow scale (or coarse-grained) dynamics which results from fast scale nonlinear dynamics.  

The SINDy algorithm has been widely applied to identify models for fluid flows~\cite{Loiseau2017jfm}, optical systems~\cite{Sorokina2016oe}, chemical reaction dynamics~\cite{Hoffmann2018arxiv},  convection in a plasma~\cite{Dam2017pf}, structural modeling~\cite{lai2019sparse}, and for model predictive control~\cite{Kaiser2018prsa}. There are also a number of theoretical extensions to the SINDy framework, including for identifying partial differential equations~\cite{Rudy2017sciadv,Schaeffer2017prsa},  and models with rational function nonlinearities~\cite{Mangan2016ieee}. It can incorporate partially known physics and constraints~\cite{Loiseau2017jfm}. The algorithm can also be reformulated to include integral terms for noisy data~\cite{Schaeffer2017pre} or handle incomplete or limited data~\cite{Tran2016arxiv,schaeffer2018extracting}. The selected modes can also be evaluated using information criteria for model selection~\cite{Mangan2017prsa}. These diverse mathematical developments provide a mature framework for broadening the applicability of the model discovery method.

Modeling multiscale physics phenomena remains difficult for modern day mathematical and computational methods.  Indeed, many complex systems exhibit diverse behaviors across multiple time and spatial scales, which poses unique challenges for characterizing and predicting their behavior. Coarse-graining techniques attempt to imbue macroscale dynamics, either in time or space, with the dominant features observed at the microscale. This can be particularly effective if the time or space scales can be clearly disambiguated in a principled way. There is a significant body of research focused on modeling multiscale systems: notably the heterogeneous multiscale modeling (HMM)  framework and equation-free methods for linking scales \cite{kevrekidis_equation-free_2003,weinan_heterognous_2003,weinan_principles_2011}. Additional work has focused on testing for the presence of multiscale dynamics so that analyzing and simulating multiscale systems is more computationally efficient \cite{froyland_computational_2014,froyland_trajectory-free_2016}. Many of the same issues that make modeling multiscale systems difficult can also present challenges for model discovery and system identification. This motivates the development of specialized methods for the discovery of coarse-grained models subject to multiple timescales.

Recently, a number of data-driven strategies for discovering nonlinear multiscale dynamical systems and their embeddings from data have been developed~\cite{Champion18}. These methods apply in two canonical cases: (i) systems for which full measurements of the governing variables are available, and (ii) systems for which there are incomplete measurements. The present framework falls into the first case as we assume here that numerical or observational trajectories of a given dynamical system have already been generated and can be used to extract data in the Poincar\'e section. For systems with incomplete observations, the Hankel alternative view of Koopman (HAVOK) method, based on time-delay embedding coordinates, can be used to obtain a linear model and Koopman invariant measurement system that nearly perfectly captures the dynamics of nonlinear quasiperiodic systems. Hence, future investigations could combine the HAVOK method with our method herein to provide a suite of mathematical strategies for reducing the data required to discover and model nonlinear multiscale systems.

This paper is laid out as follows. In Section~\ref{sec:Methods} we present all background and methods that will be used throughout this manuscript. This includes a discussion of the SINDy method in \S~\ref{subsec:SINDy} and how it can be used to discover unknown, generally nonlinear, iterative procedures. In \S~\ref{subsec:Floquet} we provide an overview of Floquet theory and its limitations as well as describe how our method can be applied to develop what we have termed nonlinear Floquet theory. Following the theoretical expositions in Section~\ref{sec:Methods}, we come to Section~\ref{sec:Applications} which provides a number of applications of our method. These applications range from simple equations to check that the methods produce equations which align with exact analytical calculations to the fully chaotic regimes of the R\"ossler differential equation. We then conclude in Section~\ref{sec:Discussion} with a discussion of our findings as well as highlight avenues for future applications and extensions of our method.

\section{Methods}\label{sec:Methods} 

This section is broken down into two subsections where we provide a description of our methods used to discover nonlinear Poincar\'e maps from data. We begin in \S~\ref{subsec:SINDy} with an overview of the SINDy method and a complete description of how it can be employed to aid in the discovery of discrete dynamical systems such as Poincar\'e maps. We conclude this section with $\S~\ref{subsec:Floquet}$ with a discussion of Floquet theory, its relation to Poincar\'e maps, and the extensions that our work provides to these theoretical dynamical systems tools.

\subsection{The SINDy Method}\label{subsec:SINDy} 

The setting of this work is to consider an iterative scheme for which given any initial condition $\xv_0 \in \mathbb{R}^d$ we may define an infinite sequence $\{\xv_0,\xv_1,\xv_2,\dots\}$ governed by the mapping
\begin{equation}\label{DDS}
	\xv_{n+1} = \Piv(\xv_n),
\end{equation}     
for some smooth function $\Piv:\mathbb{R}^d \to \mathbb{R}^d$. In our present context the function $\Piv$ will be unknown and therefore our goal is to determine accurate approximations of $\Piv$ provided we are given a finite segment of a sequence known to be governed by (\ref{DDS}). 

To arrive at accurate approximations of $\Piv$ we employ the SINDy method, as presented in \cite{SINDy}. SINDy is a sparse regression model that allows one to construct nonlinear dynamical models from data. SINDy extracts parsimonious dynamics from time-series data by taking snapshot data $\xv(t_n) \in \Rb^d$ and attempting to discover a best-fit continuous or discrete dynamical system with as few terms as possible:
\begin{equation}\label{eq:dynamical_system_x}
	\dot{\xv}(t_n) = \fv(\xv(t_n)), \quad {\rm or}\quad \xv(t_{n+1}) = \fv(\xv(t_n)), 
\end{equation}
where $\fv:\mathbb{R}^d \to \mathbb{R}^d$ is smooth. That is, the state of the system $\xv$ evolves in time $t$, with dynamics constrained by the function $\fv$. We seek a parsimonious model for the dynamics, resulting in a function $\fv$ that contains only a few active terms: it is sparse in a basis of possible functions.  This is consistent with our extensive knowledge of a diverse set of evolution equations used throughout the physical, engineering and biological sciences.   Thus, the types of functions that comprise $\fv$ are typically known from modeling experience.

SINDy frames model discovery as a sparse regression problem. If a finite sequence $\{\xv_0,\xv_1,\dots, \xv_m\}$, $m \geq 1$, assumed to satisfy (\ref{DDS}) is given, the snapshots are stacked to form data matrices $\Xv_1 = [\xv_0,\dots,\xv_{m-1}]^T$ and $\Xv_2 = [\xv_1,\dots,\xv_{m}]^T$ belonging to $\mathbb{R}^{m\times d}$. Although $\Piv$ is unknown, we can construct an extensive library of $p$ candidate functions $\Thetav(\Xv) = [\thetav_1(\Xv) \cdots \thetav_p(\Xv)] \in \Rb^{m\times p}$, where each $\thetav_j$ is a candidate model term. We assume $m\gg p$ so the number of data snapshots is larger than the number of library functions; it may be necessary to sample transients and multiple initial conditions to improve the condition number of $\Thetav(\Xv)$. The choice of basis functions typically reflects some knowledge about the system of interest: a common choice is polynomials in $\xv$ as these are elements of many canonical models. The library is used to formulate an overdetermined linear system
\[
\begin{split}
  \Xv_2 &= \Thetav(\Xv_1)\Xiv\\ &=
	\begin{bmatrix}
		\thetav_1(\xv_0) & \thetav_2(\xv_0) & \cdots & \thetav_p(\xv_0) \\
		\thetav_1(\xv_1) & \thetav_2(\xv_1) & \cdots & \thetav_p(\xv_1) \\
		\vdots & \vdots & \ddots &\vdots \\
		\thetav_1(\xv_{m-1}) & \thetav_2(\xv_{m-1}) & \cdots & \thetav_p(\xv_{m-1}) \\
	\end{bmatrix}\Xiv
\end{split}
\]
where the unknown matrix $\Xiv = [\xiv_1\ \xiv_2\ \cdots\ \xiv_n]\in \Rb^{p\times d}$ is the set of coefficients that determine the active terms from $\Thetav(\Xv)$ in the dynamics of the approximation of $\Piv$, hereby denoted $\tilde{\Piv}$. Determining an appropriate $\Xiv$ results in the form of $\tilde{\Piv}$ as a linear combination of the basis functions:
\[
	\tilde{\Piv}(\xv) = \begin{bmatrix}
		\xiv_{1,1}\thetav_1(\xv) + \dots + \xiv_{p,1}\thetav_p(\xv) \\
		\xiv_{1,2}\thetav_1(\xv) + \dots + \xiv_{p,2}\thetav_p(\xv) \\
		\vdots \\
		\xiv_{1,d}\thetav_1(\xv) + \dots + \xiv_{p,d}\thetav_p(\xv) \\
	\end{bmatrix}
\] 
which approximates the dynamics of $\Piv$. Sparsity-promoting regression is used to solve for $\Xiv$ that result in parsimonious models, ensuring that $\Xiv$, or more precisely each $\xiv_j$, is sparse and only a few columns of $\Thetav(\Xv)$ are selected. The standard SINDy approach uses a sequentially thresholded least squares algorithm to find the coefficients \cite{SINDy}, which is a proxy for $\ell_0$ optimization~\cite{zheng2019unified} and has convergence guarantees~\cite{Zhang2018arxiv}. Specifically, we promote sparsity of our solution $\Xiv$ by introducing a sparsity threshold parameter $\lambda > 0$ and defining a new unknown matrix $\tilde{\Xiv} = [\tilde{\xiv}_{i,j}] \in \mathbb{R}^{p\times d}$ with fixed $\tilde{\xiv}_{i,j} = 0$ if $|\xiv_{i,j}| \leq \lambda$ and then solving
\[
	\Xv_2 = \Thetav(\Xv_1)\tilde{\Xiv}.
\] 
We note that the matrix $\tilde{\Xiv}$ potentially has fewer degrees of freedom than $\Xiv$, depending on the value of $\lambda$. We refer the reader to the supplementary material of \cite{SINDy} which contains a MATLAB implementation of this sparsification procedure and to \url{https://faculty.washington.edu/kutz/page26/} for the full MATLAB code base.   A robust python-based SINDy package (PySindy) is also now available at {\bf GitHub/dynamicslab}. 

Critical to the success of the SINDy algorithm is the choice of initial data used to generate dynamical trajectories.  One can either use a long trajectory from a single initial condition to identify the dynamics, or alternatively use a diverse set of initial conditions to help identify the dynamics in a more efficient manner~\cite{xiu2019}. Both methods have been demonstrated to provide successful strategies, with the diversity of initial conditions ultimately requiring less data in general as it helps to identify a sparse model consistent with {\em all} the trajectories.  At this time, there are no theoretical guarantees about the quantity, diversity and quality of data required for SINDy to be successful in model discovery.  Moreover, it is highly dependent upon the problem considered.   Rudy et al.~\cite{Rudy2017sciadv} consider a number of canonical examples for which the diversity of data requirements is evaluated for model discovery in spatio-temporal systems. For the present work, we sample whenever possible from a diverse set of initial trajectories as this has been observed to lessen the data required for model discovery.

We comment on the fact that there are two major factors, or hyper-parameters, which must be considered when applying this method: the choice of functions in the library and the value of the sparsity parameter. As is pointed out in \cite{SINDy} and its subsequent investigations, these two pieces of the method are intimately related. For example, one may attempt to use only polynomial functions to approximate a purely trigonometric $\Piv(\xv)$. In this case one will not obtain a sparse relation since no finite linear combination of polynomial functions can be used to describe a trigonometric function. Choosing appropriate functions for the library can be a very difficult task and in this manuscript we do not attempt to provide general statements to guide the users choice of library functions.

Small et al.~\cite{Small2002} and Yao and Bollt~\cite{Yao2007physicad} previously formulated system identification as a similar linear inverse problem without including sparsity, resulting in models that included all terms in $\boldsymbol{\Theta}$.   In either case, an appealing aspect of this model discovery formulation is that it results in an overdetermined linear system for which many regularized solution techniques exist. Thus, it provides a computationally efficient counterpart to other model discovery frameworks~\cite{schmidt_distilling_2009}.

\subsection{Nonlinear Floquet Theory} \label{subsec:Floquet} 

In the recent multiscale SINDy architecture of Champion et al~\cite{Champion18}, sampling strategies are introduced for learning slow scale dynamics which are disambiguated from the fast scale physics.  In this work, an alternative sampling is proposed based upon Poincar\'e maps on the scale of the slow dynamics.  This provides a principled architecture for evaluating the evolution of the coarse-grained variables and discovering their slow-scale physics. 

Poincar\'e maps are intimately related to Floquet theory and the computation of the stability of periodic orbits. This is accomplished by linearizing an ordinary differential equation about a periodic orbit, resulting in a variational equation of the form
\begin{equation}\label{LinearPeriodic}
  	\dot{\bf x} = {\bf A}(t) {\bf x}
\end{equation}
where the matrix ${\bf A}$ is a square continuous function such that ${\bf A}(t+T)={\bf A}(t)$ is periodic with period $T > 0$. Then, Floquet theory dictates that any fundamental matrix $\boldsymbol{\Phi}(t)$ of this system may be decomposed as~\cite{Guckenheimer}
\begin{equation}  \label{Fundamental}
	\boldsymbol{\Phi}(t) =  \hat{\boldsymbol{\Phi}}(t) \exp( {\bf R} t)
\end{equation}
where $ \hat{\boldsymbol{\Phi}}(t) =  \hat{\boldsymbol{\Phi}}(t+T)$ is nonsingular for all $t$ and ${\bf R}$ is a constant matrix. The eigenvalues of the matrix ${\bf R}$ are the Floquet exponents which determine the stability of the periodic solutions. Specifically, if the real part of any eigenvalue is positive, solutions will blow up as $t\rightarrow \infty$. For purely imaginary eigenvalues, the solutions are bounded for all time. Floquet theory is also known as monodromy theory (Hamiltonian systems, etc.) or Bloch theory (solid state physics, etc.), where in the latter case the Floquet exponent is often referred to as the quasi-momentum. Note further that the form of the solution decomposition (\ref{Fundamental}) resembles the {\em dynamic mode decomposition} (DMD)~\cite{Kutz2018} which regresses nonlinear dynamical systems to their best-fit (least-squares) resembling the Floquet decomposition, i.e. the DMD decomposition represents the data in terms of a mode, an associated exponential time coefficient (real and imaginary parts are allowed) and a weighting of each mode.  DMD modes, unlike Floquet modes used in (\ref{Fundamental}), are not time-dependent.

A linear Poincar\'e map can be constructed from (\ref{LinearPeriodic}) using the fundamental matrix (\ref{Fundamental}) by tracking the solutions at integer multiples of the temporal period. Through this, Floquet theory is widely used for characterizing linear systems with periodic time or space coefficients. However, nonlinear systems are difficult since the decomposition form (\ref{Fundamental}) does not hold and so what is required is a method for considering the {\em nonlinear} dynamics of the evolution of the Poincar\'e map. Indeed, there is currently no general theory on the discovery or computation of nonlinear Floquet maps. This provides a natural application for our method in this manuscript. An alternative approach is to linearize the coordinate system~\cite{shai1} so that linear Floquet theory can once gain be applied~\cite{shai2}.  This alternative method instead requires computing linearizing coordinate transformations whereas our technique directly computes a dynamical model for the evolution of the map.
    
Consider a general continuous time dynamical system 
\begin{equation}\label{DS}
	\dot{\xv} = {\bf f}(\xv),	
\end{equation}
where $\xv = \xv(t) \in \mathbb{R}^{d+1}$ and ${\bf f}:\mathbb{R}^{d+1}\to \mathbb{R}^{d+1}$ is smooth. Given a trajectory of (\ref{DS}), $\xv(t)$, we may define a $d$-dimensional section $S \subset \mathbb{R}^{d+1}$ and simply track the values of $\xv(t)$ at each time it transversely intersects this section. We then define the sequence 
\[
	\{\xv_n := \xv(t_n)|\ \xv(t_n) \in S,\ t_n > t_{n-1}\} 
\] 
and consider the iterative scheme which takes $\xv_n \mapsto \xv_{n+1}$. This iterative scheme defines a Poincar\'e map that can be used to determine the existence and stability of periodic orbits and invariant manifolds as well as providing a dimensional reduction of chaotic dynamics. We see that in this more general context such a Poincar\'e map generalizes Floquet theory to a range of dynamical systems which are neither linear nor periodic in the independent variable.  

In practice a mapping which takes $\xv_n \mapsto \xv_{n+1}$ can rarely be found explicitly, especially away from the neighborhood of a known solution such as a periodic orbit. Presently we are left to inspect simplified models which can explain certain nonlinear phenomena, but cannot forecast dynamics of specific dynamical systems. This problem can now potentially be overcome by using the SINDy method for maps. That is, we may numerically integrate (\ref{DS}) and generate a sequence of training data by constructing a sequence of iterates which lie in the prescribed section at successive times. In this case we arrive at an unknown mapping (\ref{DDS}) which we wish to discover via the SINDy method. In the following section we will provide a number of examples for which we can apply this method to dynamical systems which exhibit both simple and complex dynamical behavior.   

Importantly, the training data is used to derive a parsimonious representation of the nonlinear dynamics which has the strongest potential for generalization.  Indeed, the fewer the number of parameters representing the dynamics and data, the higher the probability of producing a generalizable model.  This is consistent with long standing ideas of the Pareto front, parsimony and interpretability~\cite{BK2019}.

\section{Applications of the Method}\label{sec:Applications} 

In this section we go through some important examples of how our methods can be applied. Throughout we will use a library containing only constant functions and polynomial terms up to degree five. More precisely, we take 
\[
	\Thetav(\xv) = [1, \mathbb{X}, \mathbb{X}^{P_2}, \mathbb{X}^{P_3}, \mathbb{X}^{P_4}, \mathbb{X}^{P_5}] .
\]
Here we have used the notation $\mathbb{X}^{P_i}$ to denote the homogeneous $i^{th}$ degree polynomials in $\xv \in \mathbb{R}^d$. For example, the linear terms are given by
\[
	\mathbb{X} = [\xv_1,\dots,\xv_d],
\]
and the quadratic polynomials are given by
\[
	\mathbb{X}^{P_2} = [\xv_1^2, \xv_1\xv_2, \cdots, \xv_2^2, \cdots, \xv_d^2].
\]  
The specific value of $d\geq 1$ will vary between examples and will be clear from the context.

\subsection{RC Circuit} 

We consider the simple example presented in \cite[Example 8.7.2]{Strogatz} for which a Poincar\'e map can be determined analytically. An RC circuit is modeled in dimensionless form by the non-autonomous ordinary differential equation 
\begin{equation}\label{RC}
	\dot{x} = \sin(2\pi t) - x.
\end{equation}
Since the forcing term is $1$-periodic, a natural Poincar\'e map is obtained by mapping $x(0)$ to $x(1)$. Solving (\ref{RC}) exactly leads to the Poincar\'e mapping given by 
\begin{equation}\label{RC_Map}
	\Piv(x) = \mathrm{e}^{-1}x - \frac{2\pi(\mathrm{e} - 1)}{\mathrm{e}(1 + 4\pi^2)} \approx 0.36788x - 0.09812,
\end{equation}
which is linear and exhibits a globally attracting fixed point at $x_* = -2\pi/(1 + 4\pi^2)$, corresponding to a globally attracting periodic orbit of (\ref{RC}). 

To apply our method we obtain five solutions to the differential equation (\ref{RC}) with random initial conditions at $t = 0$ and iterating forward to approximately $t = 10^5$ to guarantee convergence to the unique attractor. This results in $5\times 10^5$ data points to seed the SINDy algorithm. Then, tracking the values of these solutions at integer values of the independent variable $t$ we are able to generate the data required to apply our method. Our results are summarized in Table~\ref{ta:rc}.

\begin{center}
\begin{table}[t]
\begin{tabular}{|c|c|}
\hline
Sparsity Parameter ($\lambda$) & Poincar\'e map ($\tilde{\Piv}(x)$)  \\ \hline
0.1 & $-0.11349x^2 + 0.80631x $ \\ 
0.05 & $0.36326x -0.08599$ \\ 
0.01 & $0.36326x- 0.08599$ \\
0.005 & $0.36326x-0.08599$ \\
0.001 & $-0.00318x^5 + 0.03473x^4 -0.12645x^3 + 0.17043x^2 + 0.30203x -0.09984$ \\ \hline
\end{tabular}
\caption{\label{ta:rc}  Summary of results on the RC  differential equation (\ref{RC}).}
\end{table}
\end{center}

We can see from this table that the value of $\lambda$ is very influential on the mapping which we obtain via our method. We see that it is possible to lose the constant term by taking $\lambda$ too large, thus providing a quadratic mapping $\tilde{\Piv}$ which bears little resemblance to the true map (\ref{RC_Map}). Similarly, taking $\lambda$ small introduces unnecessary terms that increase the polynomial order of the obtained Poincar\'e map, but add little to the qualitative dynamics of the mapping. Hence, when applying this method one should be careful since the specific tuning of the parameter $\lambda$ can be problem dependent.  The process of tuning $\lambda$ can potentially be automated by considering models near the Pareto frontier along with information criteria for model selection~\cite{BK2019}.

\subsection{Hopf Normal Form} 

We now consider a simple autonomous system whose Poincar\'e mapping can again be computed explicitly. We consider the normal form of a supercritical Hopf bifurcation truncated at cubic order given by
\begin{equation}\label{Hopf}
	\begin{split}
		\dot{x} &= x - \omega y - x(x^2 + y^2), \\
		\dot{y} &= \omega x + y - y(x^2 + y^2), \\ 
	\end{split}
\end{equation} 
where $\omega \in \mathbb{R}\setminus\{0\}$ is taken to be a constant. Here we take the positive $x$-axis to be a Poincar\'e section, for which in this case we the associated Poincar\'e mapping is given by
\begin{equation}\label{Hopf_Map}
	\Piv(x) = [1 + \mathrm{e}^{-4\pi/\omega}(x^{-2} - 1)]^{-\frac{1}{2}}.
\end{equation}
In this case we have exactly two fixed points: a trivial point $x = 0$ corresponding the equilibrium at the origin, and a fixed point at $x = 1$ corresponding to the unique stable limit cycle of (\ref{Hopf}). Moreover, any positive initial condition provided to the map (\ref{Hopf_Map}) iterates monotonically to the stable fixed point at $x = 1$, thus indicating the global stability of the limit cycle. 

\begin{figure} 
\center
\includegraphics[height=0.4\textwidth]{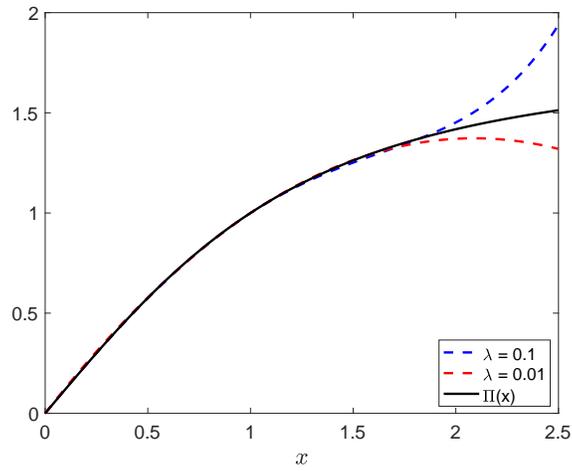}
\caption{Poincar\'e maps for the Hopf normal form (\ref{Hopf}). In black is the true map (\ref{Hopf_Map}) and its approximations $\Pi(x)$ with (blue) $\lambda = 0.1$ and (red) $\lambda = 0.01$. All three maps have fixed points at $x = 0$ and near $x = 1$.}
\label{fig:Hopf_Comparison}
\end{figure}

Although this Poincar\'e map is known explicitly, our methods will be unlikely to discover it since one would require the exact right-hand-side of (\ref{Hopf_Map}) to be in the library. Nonetheless, we may apply our method to obtain Poincar\'e mappings that qualitatively capture the dynamics of the map (\ref{Hopf_Map}). To illustrate this, we again seed our mapping with data from five trajectories obtained with random initial conditions and iterate forward for at least $300$ full periods to generate the training data. Recall that we are using a library containing the constant functions and polynomials in $x$ up to the fifth degree, and so using a sparsity parameter value of $\lambda = 0.1$ results in the quadratic approximation of (\ref{Hopf_Map}), given by
\[
	\tilde{\Piv}(x) = 1.3143x -0.31465x^2,
\]
which has fixed points at $x = 0$ (unstable) and $x = 0.9990$ (stable). Although this mapping is not an exact representation of (\ref{Hopf_Map}), it very nicely captures the qualitative dynamics of the true Poincar\'e mapping in that initial conditions taken between the roots of the polynomial will monotonically converge to the stable fixed point after at most one iteration. One may decrease the value of the sparsity parameter to obtain increasingly better approximations which use more and more polynomial values. For example, repeating the above with $\lambda = 0.01$ gives a quintic approximation of (\ref{Hopf_Map}), given by  
\[
	\tilde{\Piv}(x) = 1.2392x -0.107x^2 -0.15203x^3  + 0.019308x^5.
\]
Here we see that we have increased the number of fixed points, but again we find one that is stable and positive at $x = 0.9984$. Furthermore, the quintic map now provides monotonic convergence of any positive initial condition below the rightmost unstable equilibrium to the stable equilibrium. Therefore, we see that even with a limited library of polynomial terms our method is able to provide a qualitative approximation of the true Poincar\'e map. A comparison of the two approximate maps against the true map (\ref{Hopf_Map}) is provided in Figure~\ref{fig:Hopf_Comparison}. 

\subsection{Singular Perturbations and Averaging} 

\begin{figure} 
\center
\includegraphics[height=0.4\textwidth]{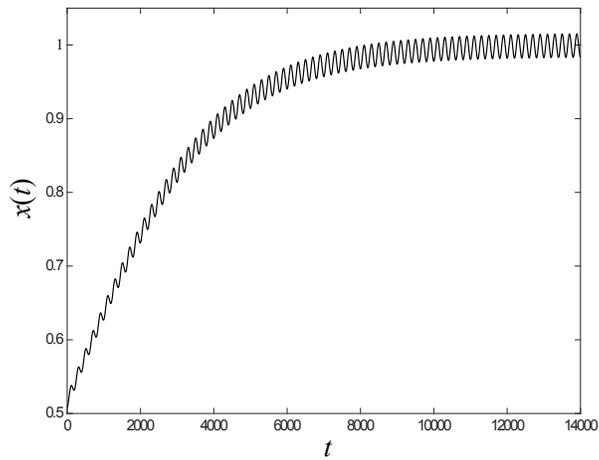}
\caption{Numerical integration of (\ref{Logistic}) with $\varepsilon = 0.1$ showing trajectories converging to a stable limit cycle.}
\label{fig:Logistic_Sol}
\end{figure} 

The methods of this work can also be applied to differential equations which exhibit multiple timescales. To illustrate this point we focus on a singularly perturbed logistic function with seasonal fluctuations. The model is given by
\begin{equation}\label{Logistic}
	\dot{x} = \varepsilon x(1 + \sin(2\pi t) - x),
\end{equation}
where $\varepsilon > 0$ is taken to be a small parameter in the system. Since (\ref{Logistic}) is $1$-periodic in $t$, one can apply the method of averaging (see \cite[Section 4.1]{Guckenheimer}) to find that the system exhibits an unstable equilibrium at $x = 0$ and a stable periodic orbit which attracts all solutions with positive initial condition. This scenario is exemplified in Figure~\ref{fig:Logistic_Sol} where we present a numerical integration of the equation (\ref{Logistic}) with $\varepsilon = 0.1$.   

As in the first example, our Poincar\'e map is obtained by tracking the value of the numerically integrated solution at integer values of $t$ since the differential equation (\ref{Logistic}) is $1$-periodic. We apply our method to (\ref{Logistic}) with $\varepsilon = 0.1$ by simulating the system for $t \in [0,15000]$ (giving $15000$ data points for the Poincar\'e section) to show that the resulting Poincar\'e map informs the global structure described above by the method of averaging. As in the previous examples, we seed our numerical integration with five randomly chosen initial conditions. A sparsity parameter of $\lambda = 0.01$ results in the approximated Poincar\'e map 
\[
	\tilde{\Piv}(x) = 1.0903x - 0.091639x^2. 
\] 
We can see that our map $\tilde{\Piv}$ has exactly two fixed points, $x = 0$ and $x = 0.9855$, with the former being unstable and the latter stable.

\subsection{Driven Brusselator} 

\begin{figure} 
\center
\includegraphics[width=0.45\textwidth]{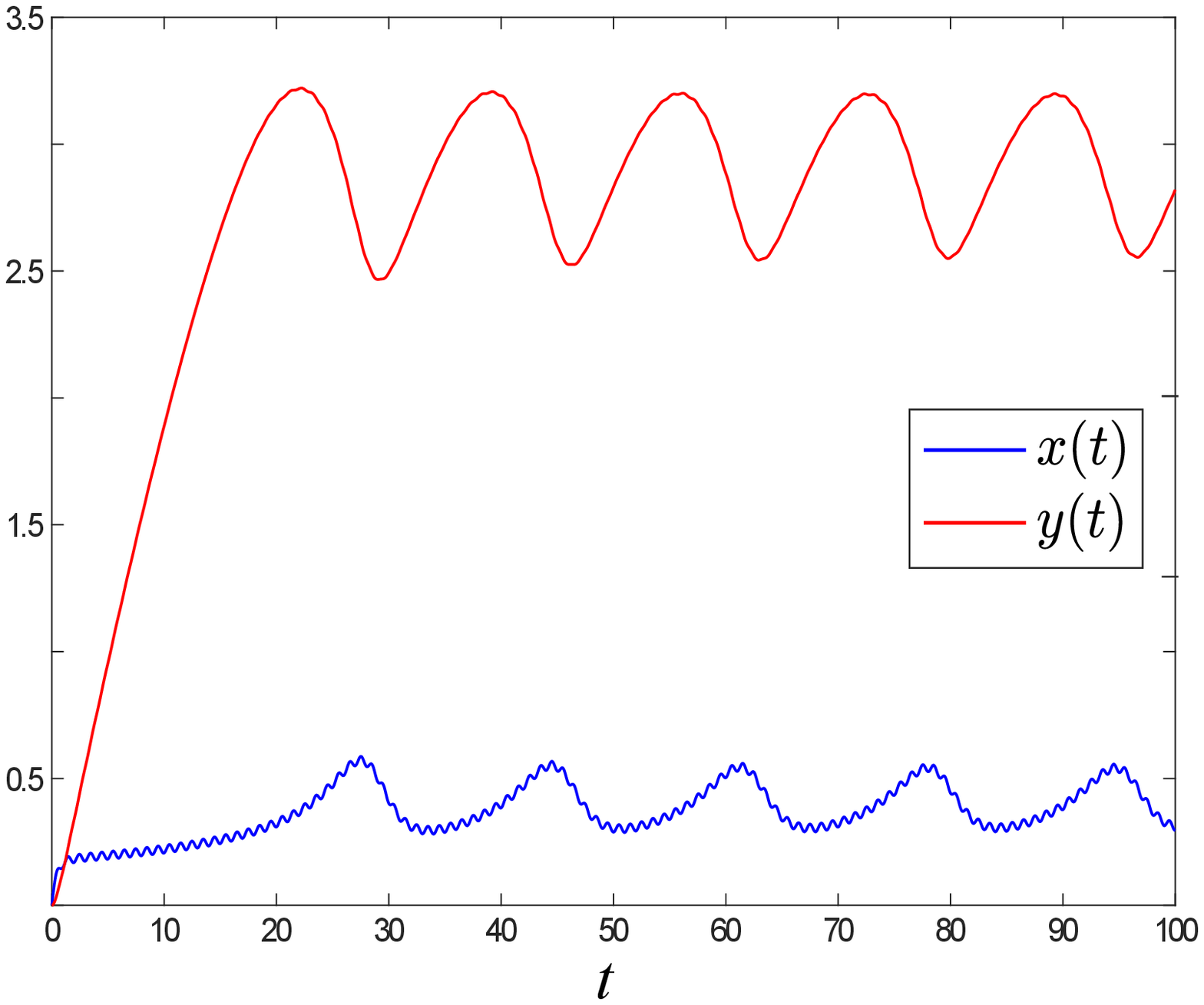} \quad
\includegraphics[width=0.45\textwidth]{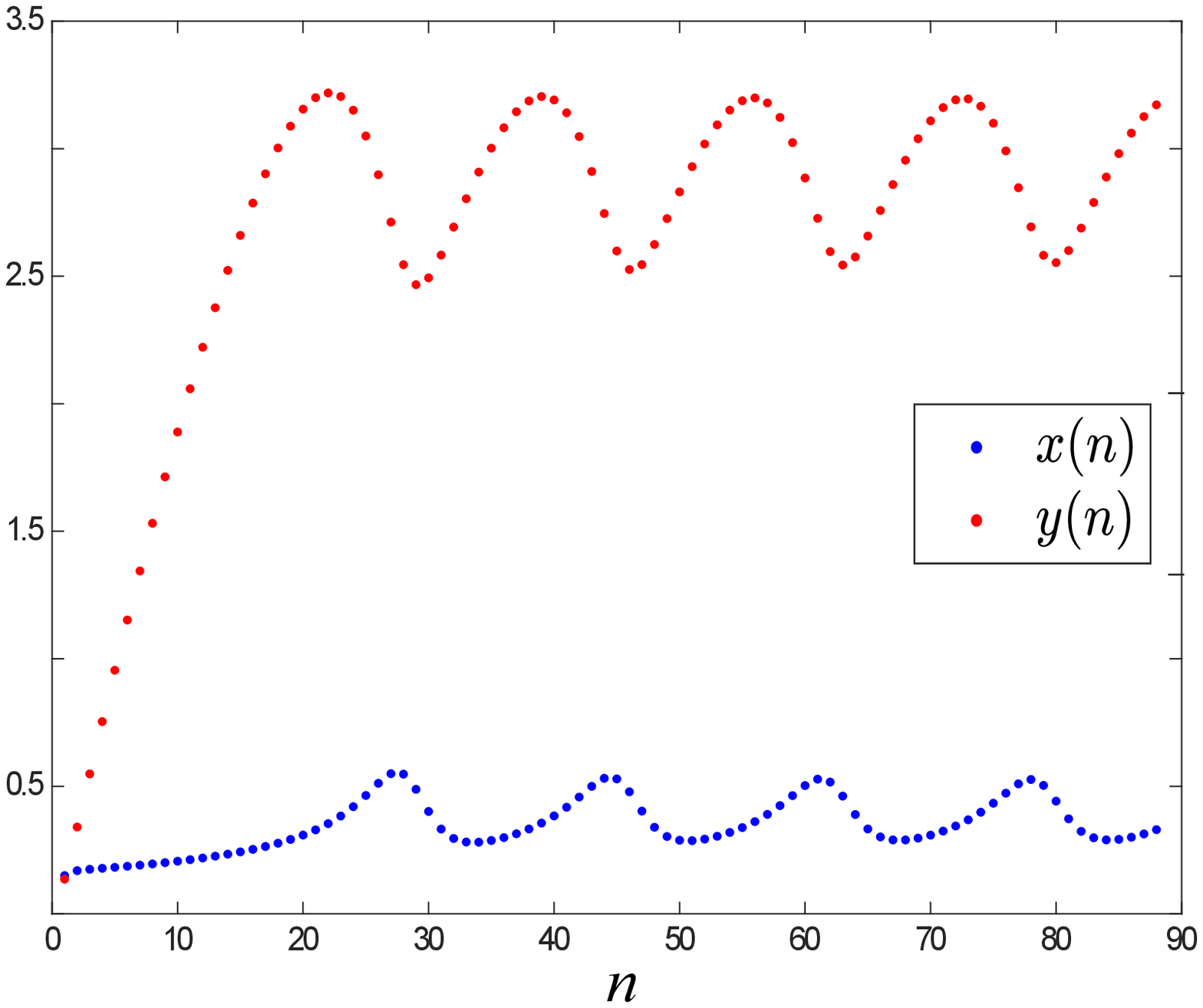} 
\caption{Numerical integration (left) and Poincar\'e section data (right) of the forced Brusselator system (\ref{Brusselator}) with initial condition $(x(0),y(0)) = (0,0)$. The parameters are taken to be $(a,b,\alpha) = (0.4,1.2,0.1)$.}
\label{fig:Brusselator_Sol}
\end{figure} 

Up to this point we have only focussed on constructing Poincar\'e maps for systems whose training data reveals only the existence of fixed points. With this example we apply our methods to illustrate that systems which exhibit quasi-periodic attractors can be constructed from training data. We will focus on the two-dimensional non-autonomous driven Brusselator, given by
\begin{equation}\label{Brusselator}
	\begin{split}
		\dot{x} &= a + \alpha\sin(2\pi t) - (b + 1)x + x^2y \\
		\dot{y} &= bx - x^2y,
	\end{split}
\end{equation}   
where $a,b,\alpha \in \mathbb{R}$ are parameters. This system models a theoretical chemical reaction with a periodically supplied reactant, represented by the $\sin(2\pi t)$ term. When $\alpha = 0$ it is well-known that there exists a connected, open set of parameter values for which the system (\ref{Brusselator}) falls into a periodic orbit \cite[]{Strogatz}. Hence, in this parameter range one would expect that slowly introducing $\alpha > 0$ will cause the system to fall into a quasi-periodic attractor. This is exemplified in Figure~\ref{fig:Brusselator_Sol} where we have provided a numerical integration of (\ref{Brusselator}) with initial condition $(x(0),y(0)) = (0,0)$ and parameter values $(a,b,\alpha) = (0.4,1.2,0.1)$. On the right in Figure~\ref{fig:Brusselator_Sol} we provide the data used to generate the Poincar\'e map which is obtained by tracking the value of the solutions $(x(t),y(t))$ at $t = n \in \mathbb{Z}$. The training data is taken so that approximately 4 full periods can be observed on the attractor.  

\begin{figure} 
\center
\includegraphics[width=0.45\textwidth]{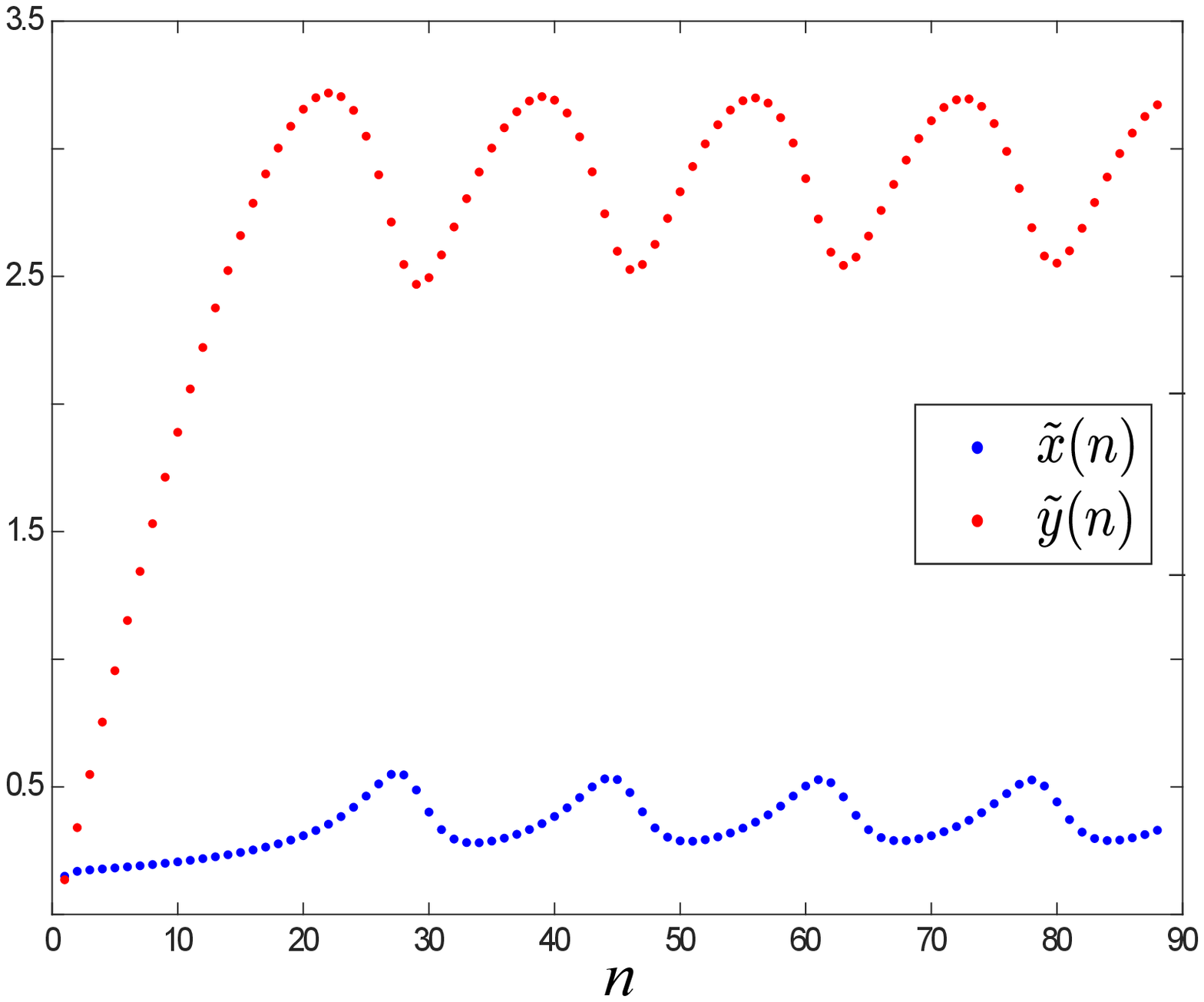} \quad
\includegraphics[width=0.45\textwidth]{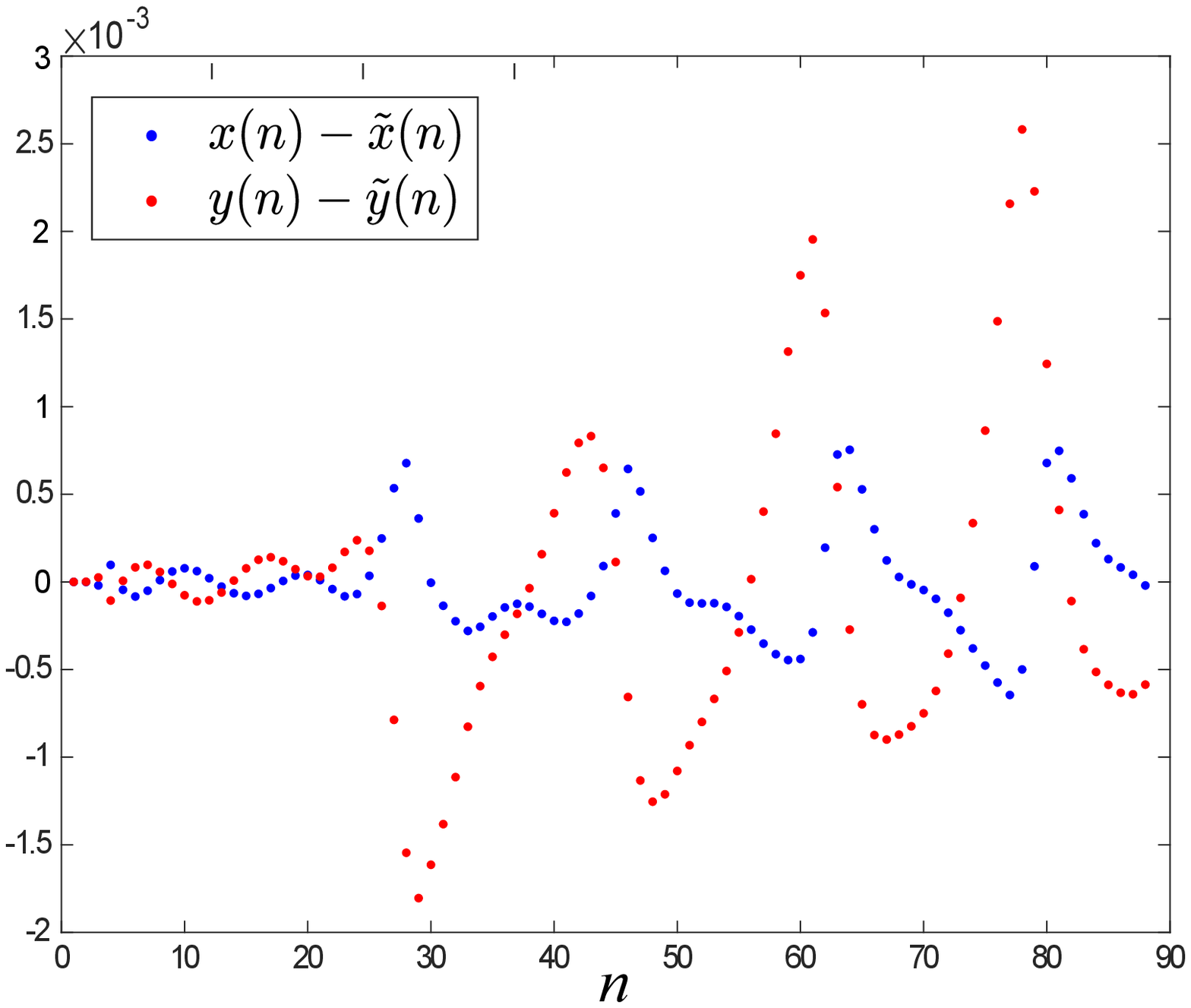} 
\caption{Trajectories of the approximate Poincar\'e map (left), denoted $(\tilde{x}(n),\tilde{y}(n))$, and their error against the original training data (right). Here we have used the data from the right of Figure~\ref{fig:Brusselator_Sol} to apply our method and obtain an approximate Poincar\'e mapping for the system (\ref{Brusselator}).}
\label{fig:Brusselator_Psec}
\end{figure} 

\begin{figure} 
\center
\includegraphics[width=0.45\textwidth]{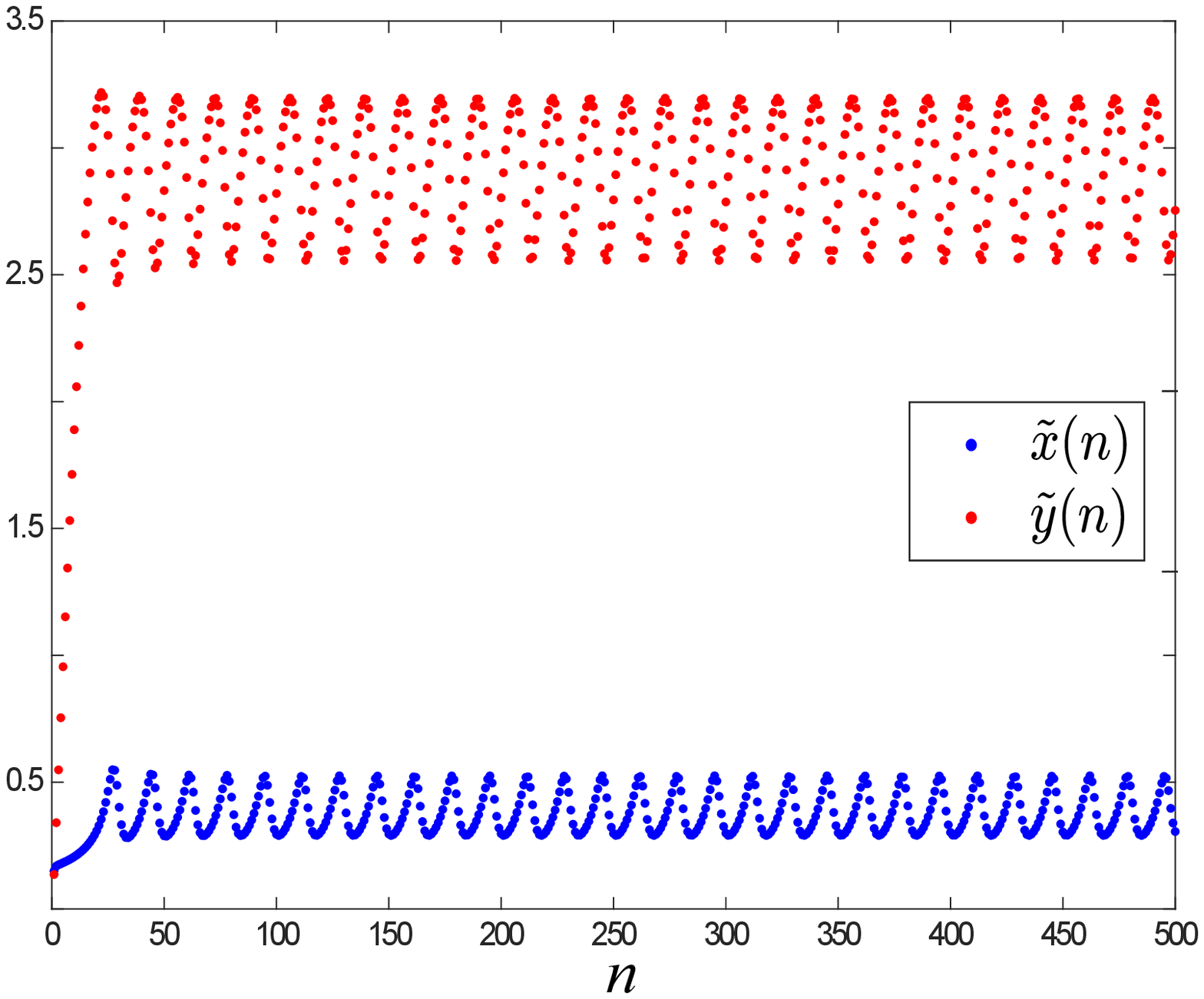} \quad
\includegraphics[width=0.45\textwidth]{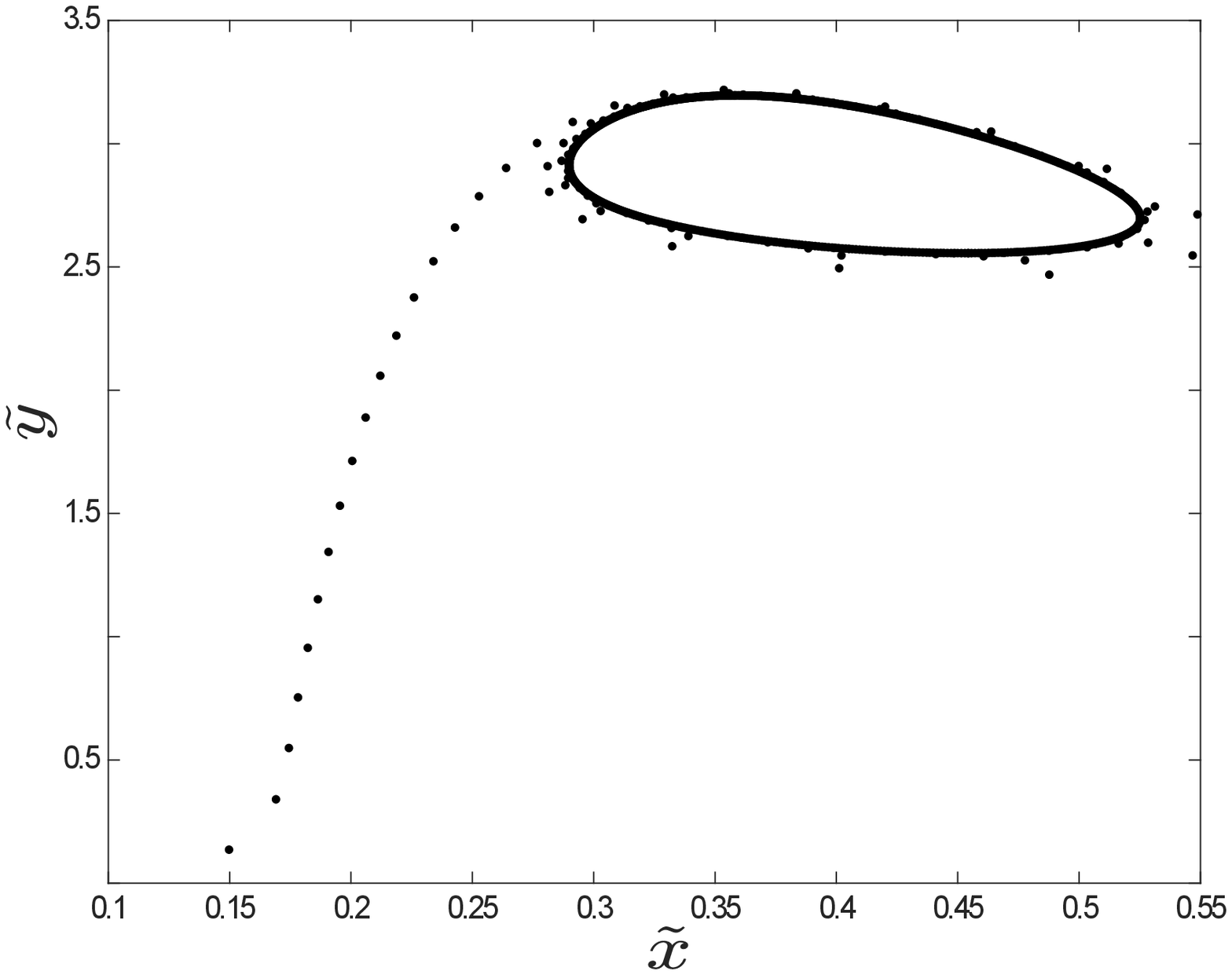} 
\caption{Longer time trajectories of the approximate Poincar\'e map associated to the forced Brusselator (\ref{Brusselator}). Presented are the time series (left) and iterations in the $(\tilde{x},\tilde{y})$-plane (right), both showing convergence to a quasi-periodic attractor.}
\label{fig:Brusselator_Psec2}
\end{figure}

On the left of Figure~\ref{fig:Brusselator_Psec} we present a plot of solutions of the approximated Poincar\'e map, denoted $(\tilde{x}(n),\tilde{y}(n))$, generated using our method with the sparsity parameter $\lambda = 0.01$. Our initial conditions are taken to be the same as the initial conditions of the input data, and hence we can see how well our approximate mapping follows the true Poincar\'e data by comparing with Figure~\ref{fig:Brusselator_Sol}. In particular, on the right of Figure~\ref{fig:Brusselator_Sol} we provide a plot of the error at each iteration of the true data versus the approximate data generated with our mapping. As one can see, the error remains on the order of $10^{-3}$, one full order of magnitude below the sparsity parameter $\lambda$. Furthermore, the error in $x(n) - \tilde{x}(n)$ appears to remain confined to a region near zero, whereas the error in $y(n) - \tilde{y}(n)$ appears to be growing as the number of iterations increases. This growth in the error can be attributed to the approximate Poincar\'e mapping not completely capturing the microscopic dynamics on the attractor, and hence $\tilde{y}(n)$ appears to fall behind or jump ahead of the training data $y(n)$. In Figure~\ref{fig:Brusselator_Psec2} we provide time series for $(\tilde{x}(n),\tilde{y}(n))$ with $n$ running from $0$ to $500$. Here it is apparent that the dynamics of the approximated Poincar\'e map does indeed rapidly converge to a quasi-periodic attractor, as is expected from the preceding discussion. This boundedness of $\tilde{y}(n)$ further implies that the error in $y(n) - \tilde{y}(n)$ is bounded and likely quasi-periodic as well.    

\subsection{R\"ossler System} 

We now focus on determining approximate Poincar\'e maps for a system which can exhibit chaotic behavior. We use the R\"ossler system as our motivating example, given by
\begin{equation}\label{Rossler}
	\begin{split}
		\dot{x} &= -y - z \\
		\dot{y} &= x + ay \\
		\dot{z} &= b + z(x-c),
	\end{split}
\end{equation}
where $a,b,c \in \mathbb{R}$ are taken to be parameters. Throughout this section we will fix $a = b = 0.1$ and vary $c$. The R\"ossler equations have the property that a solution crosses the plane $x = 0$ when $z = 0$. Therefore we take the $x = 0$ plane as our Poincar\'e section so that we are able to restrict ourselves to one-dimensional Poincar\'e maps which track the value of $y$ as $x$ crosses $0$ with positive slope.   

\begin{figure} 
\center
\includegraphics[width=0.45\textwidth]{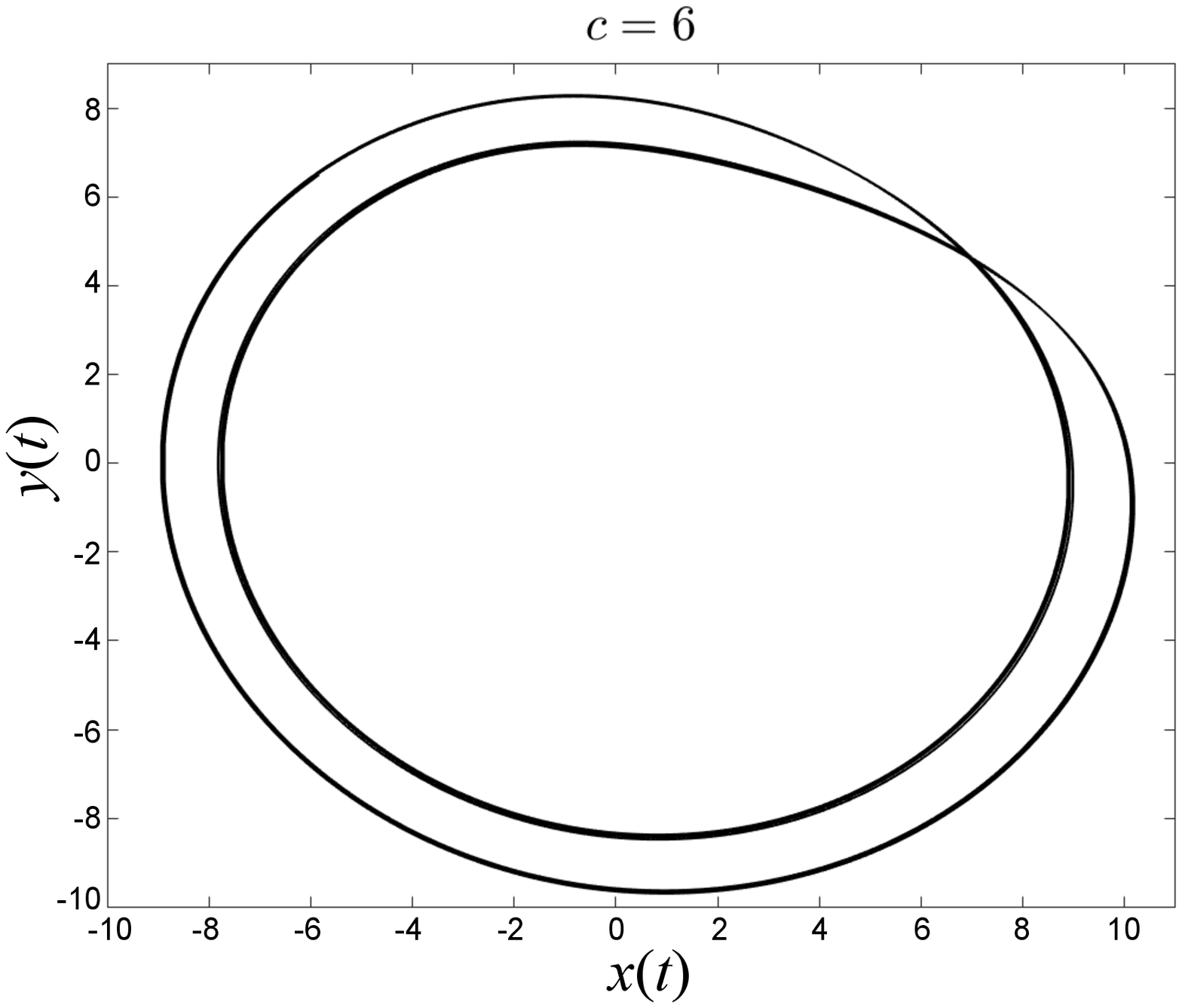} \quad
\includegraphics[width=0.45\textwidth]{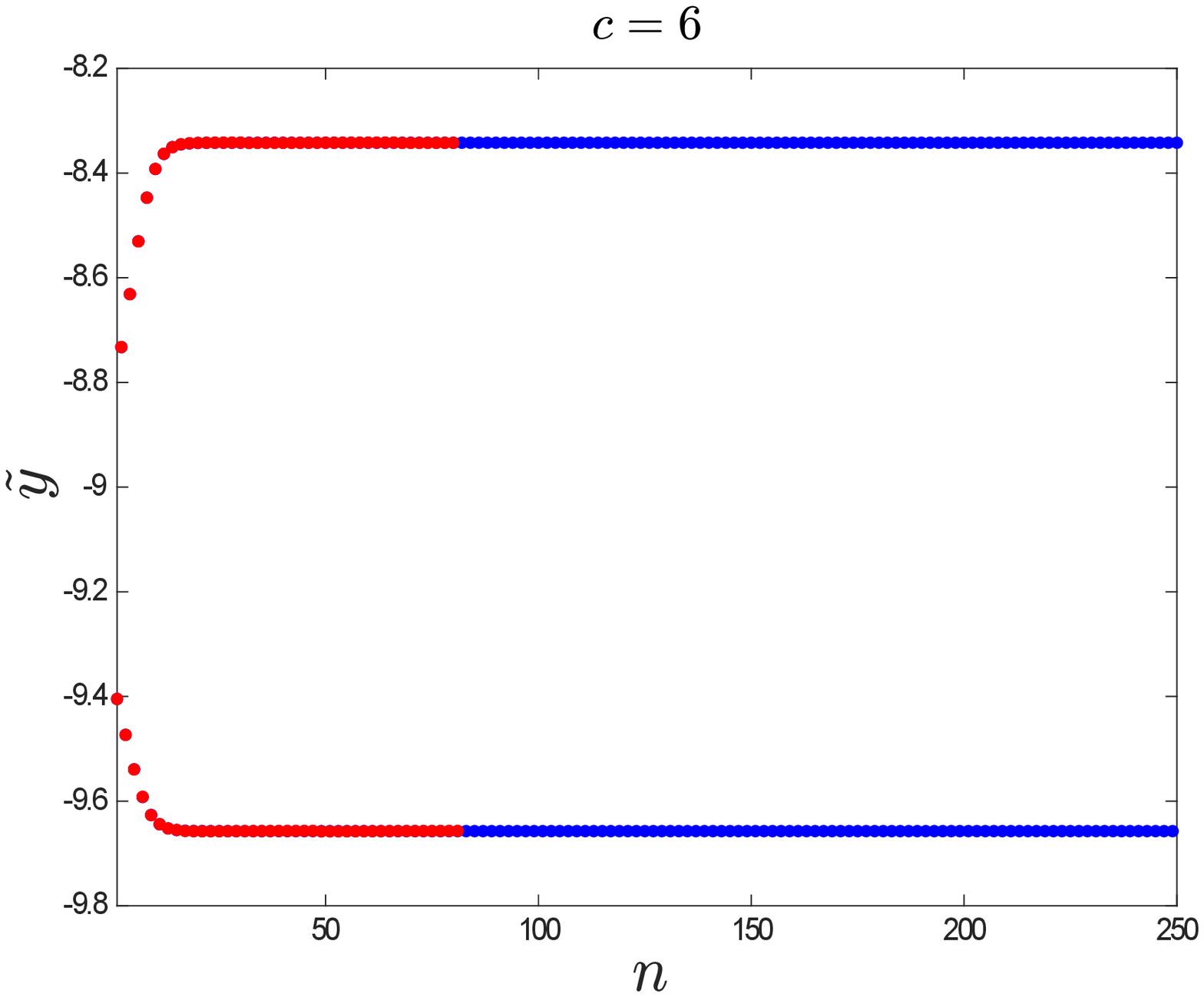}  \\
\includegraphics[width=0.45\textwidth]{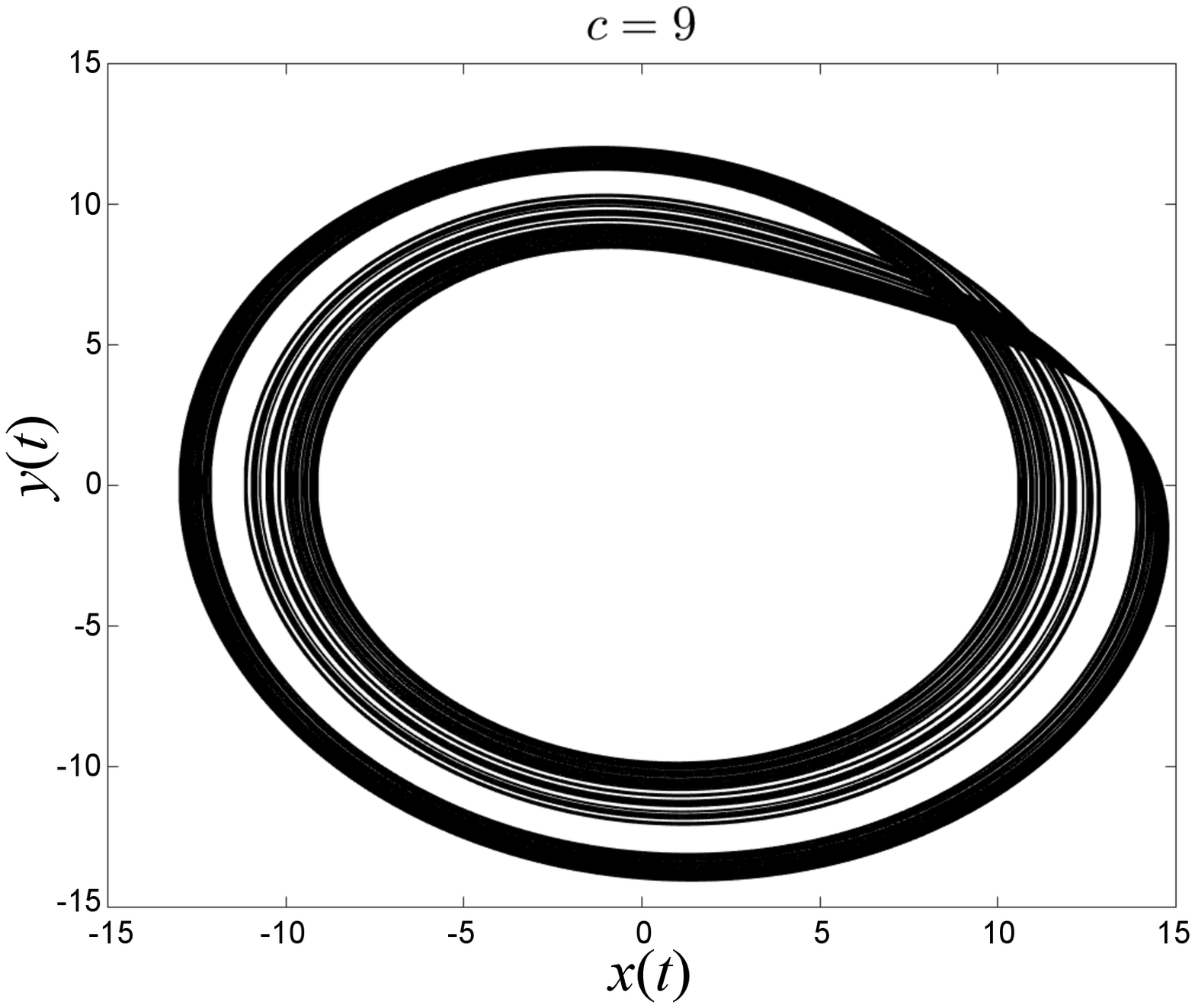} \quad
\includegraphics[width=0.45\textwidth]{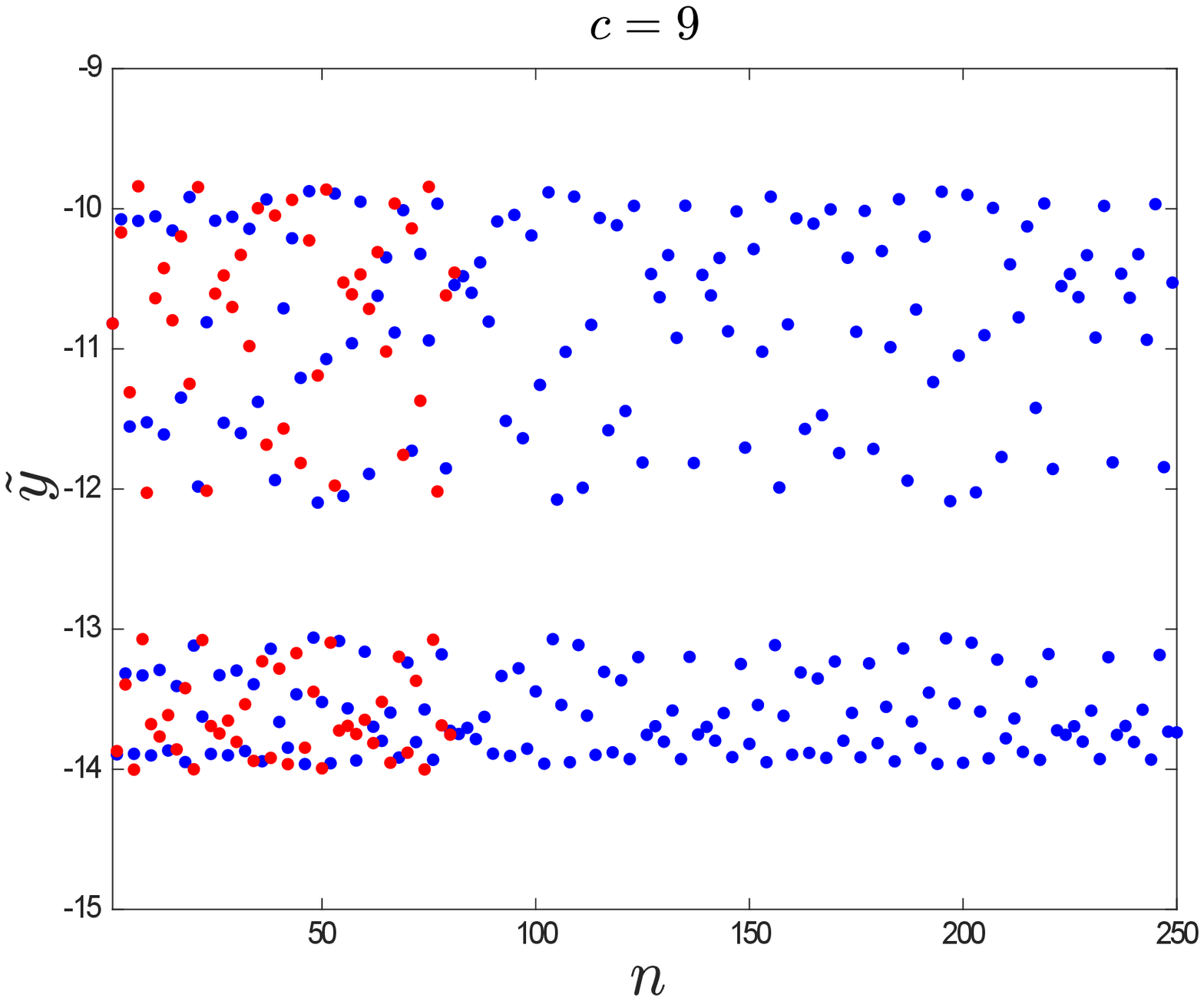} \\
\includegraphics[width=0.45\textwidth]{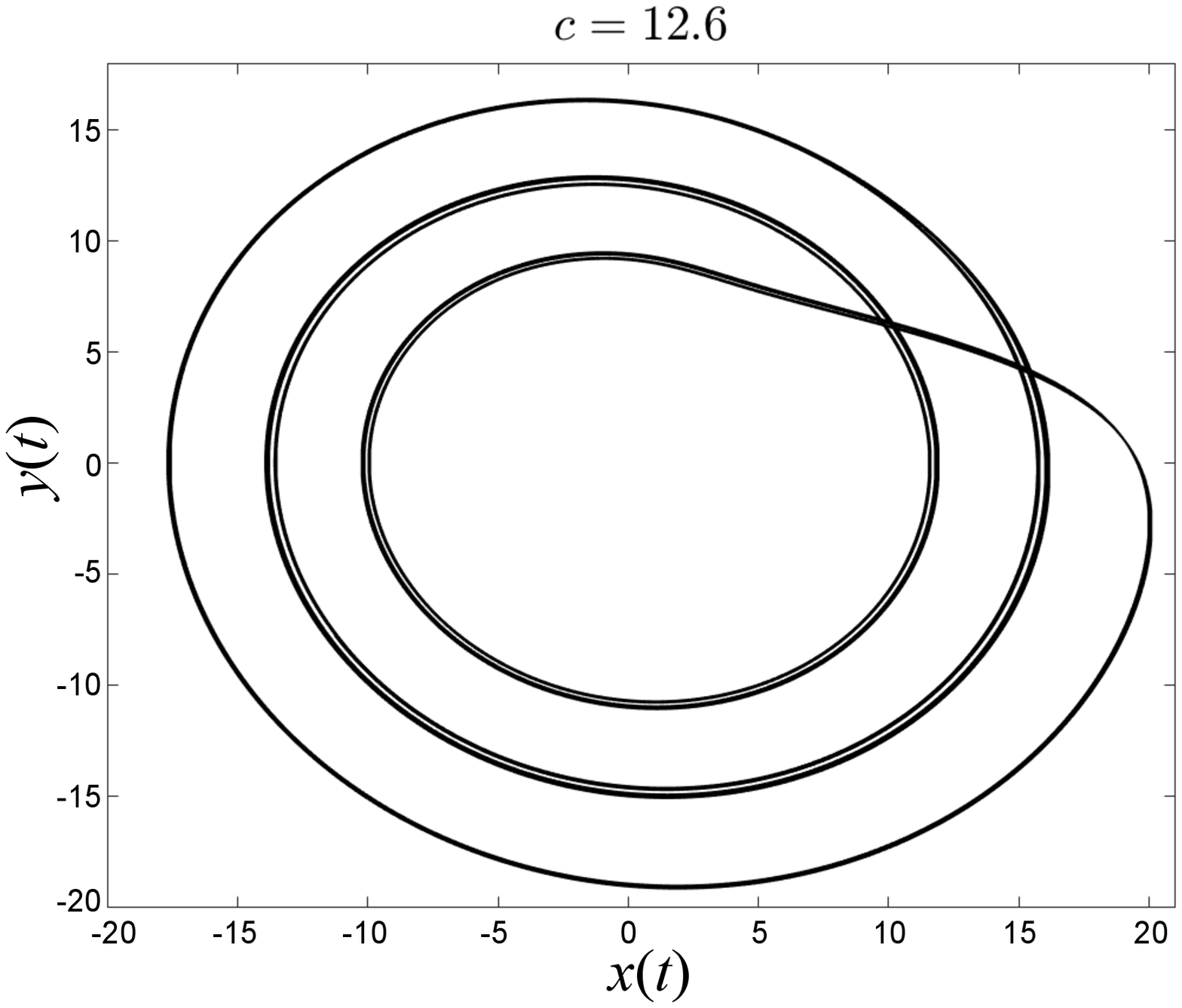} \quad
\includegraphics[width=0.45\textwidth]{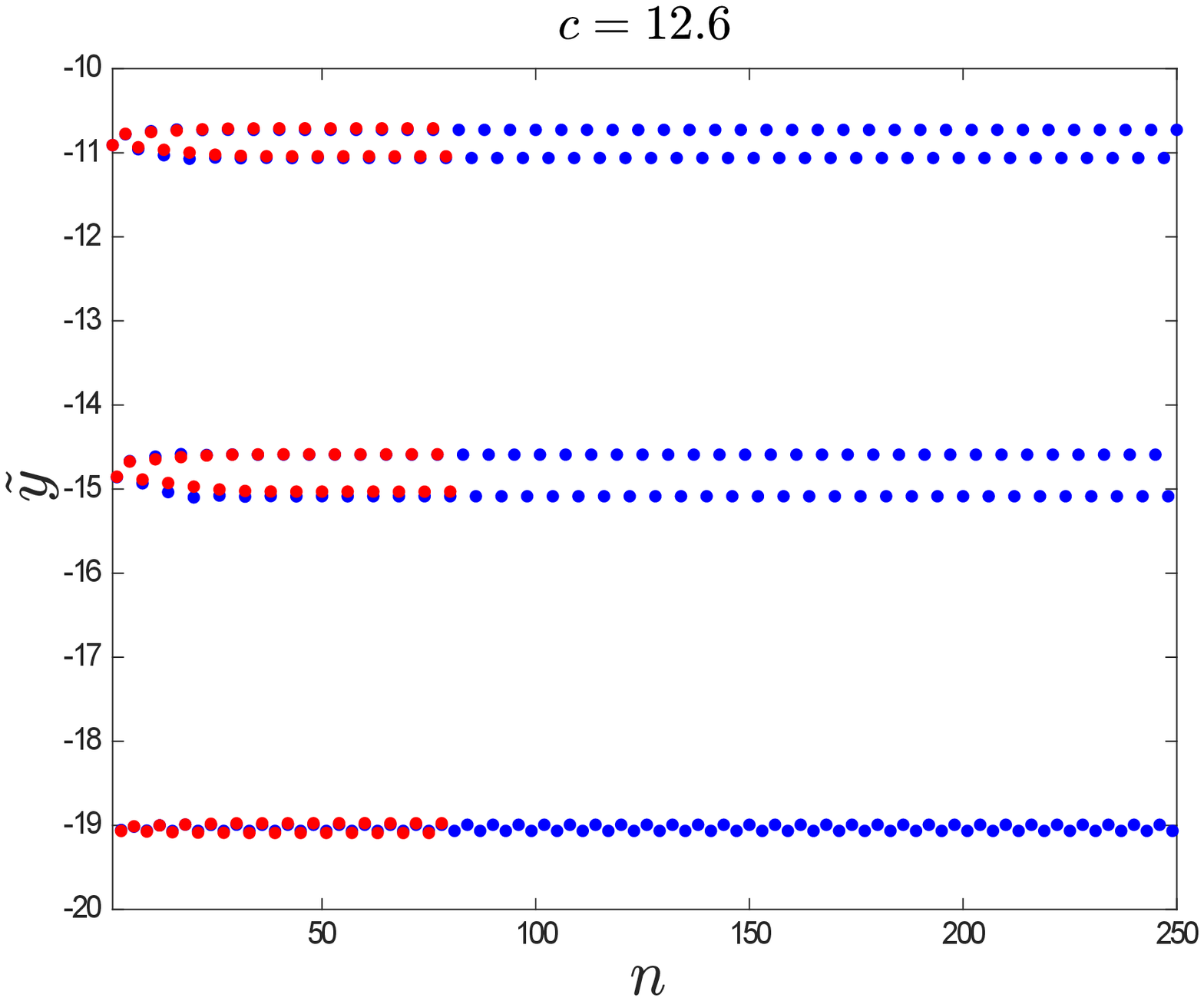}  \\
\includegraphics[width=0.45\textwidth]{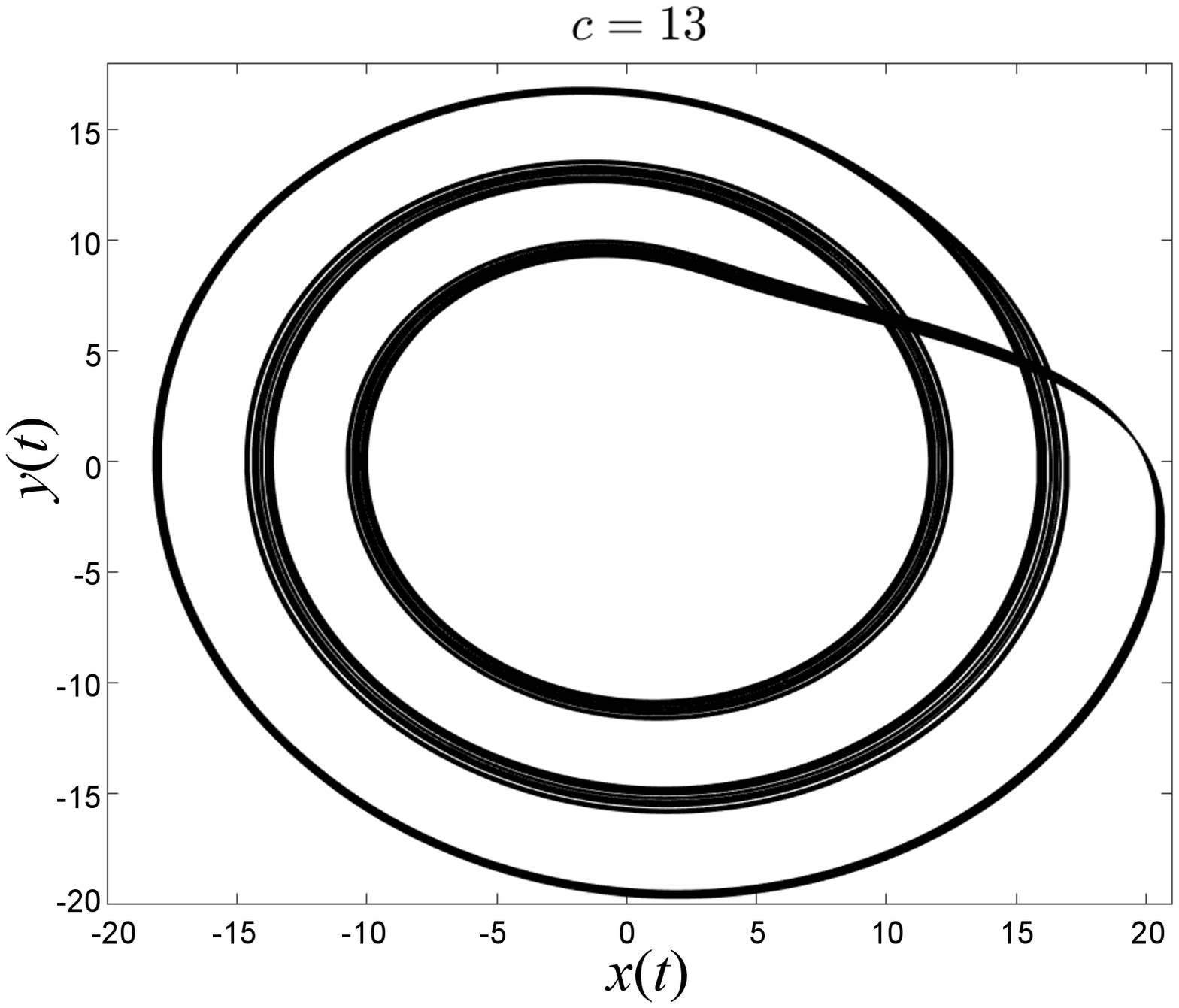} \quad
\includegraphics[width=0.45\textwidth]{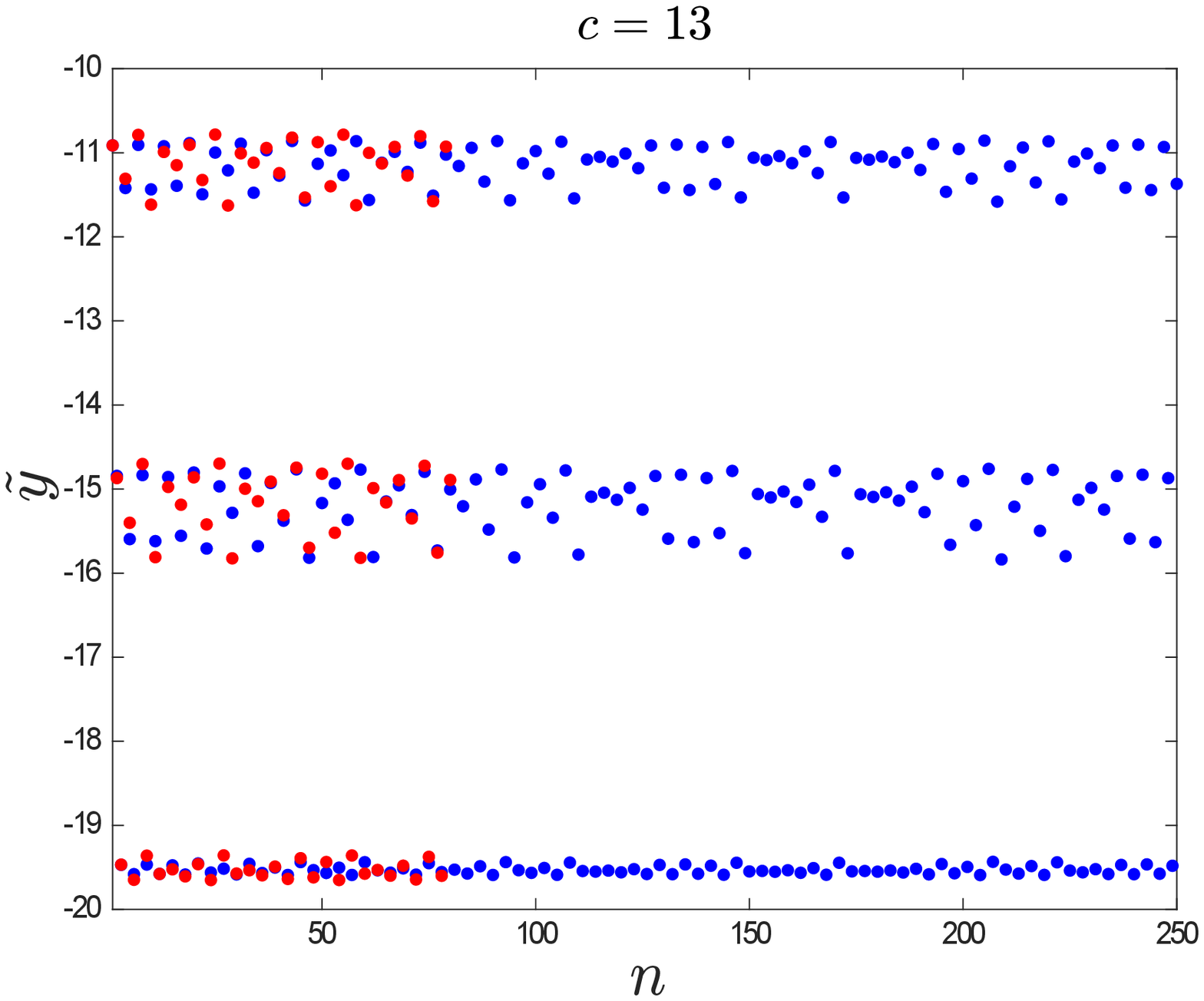} 
\caption{The R\"ossler attractor in the $(x,y)$-plane (left) and its simulated Poincar\'e map (right) for varying parameter values $c$ and a fixed sparsity parameter $\lambda = 0.01$. On the right we provide the training data in the Poincar\'e section (red) along with long time simulated dynamics of the Poincar\'e mapping using the methods of this manuscript (blue). Parameter values are given from top to bottom by $c = 6$ (2-cycle), $9$ (chaos), $12.6$ (6-cycle), and $13$ (chaos), respectively.}
\label{fig:Rossler}
\end{figure} 

Our discussion begins with describing the dynamics of (\ref{Rossler}) at four particular parameter values of $c$ which will be the focus of this example. At $c = 6$ the attractor intersects the $x = 0$ plane at exactly two points, giving an attracting 2-cycle in the Poincar\'e map. At $c = 9$ the attractor is chaotic and exhibits two distinct `bands' which are apparent in the Poincar\'e section. At $c = 12.6$ the attractor intersections the Poincar\'e section at exactly six points, giving an attracting 6-cycle in the Poincar\'e map. Finally, at $c = 13$ the attractor is again chaotic, but now exhibits three `bands' which are again apparent in the Poincar\'e section. In Figure~\ref{fig:Rossler} we provide visualizations of the attractor in the $(x,y)$-plane at these parameter values on the left and the numerically obtained training data from the Poincar\'e section on the right in red.   

On the right of Figure~\ref{fig:Rossler} we further provide long-time dynamics of the sparse Poincar\'e mapping (given in blue) obtained by applying our methods with a sparsity parameter of $\lambda = 0.01$ and approximately 75 periods of data on the attractor. We can see that in the cases of $c \in \{6,12.6\}$ where the attractor gives a cycle in the Poincar\'e section our method produces a mapping with a nearly identical attractor. In the cases when $c \in \{9,13\}$ and the attractor is chaotic, we do not expect that our mapping would be able to completely reproduce the chaotic dynamics, but we do point out that it almost perfectly captures the banding phenomenon which the R\"ossler attractor is known for. Furthermore, the defining characteristic of chaos, {\em sensitivity to initial data}, implies that one would be met with a lot of difficultly completely reproducing the dynamics of the training data. Hence, the fact that our sparse mapping is able to completely replicate the banded chaos is much more important that completely reproducing the dynamics of the training data.

\begin{figure} 
\center
\includegraphics[width=0.45\textwidth]{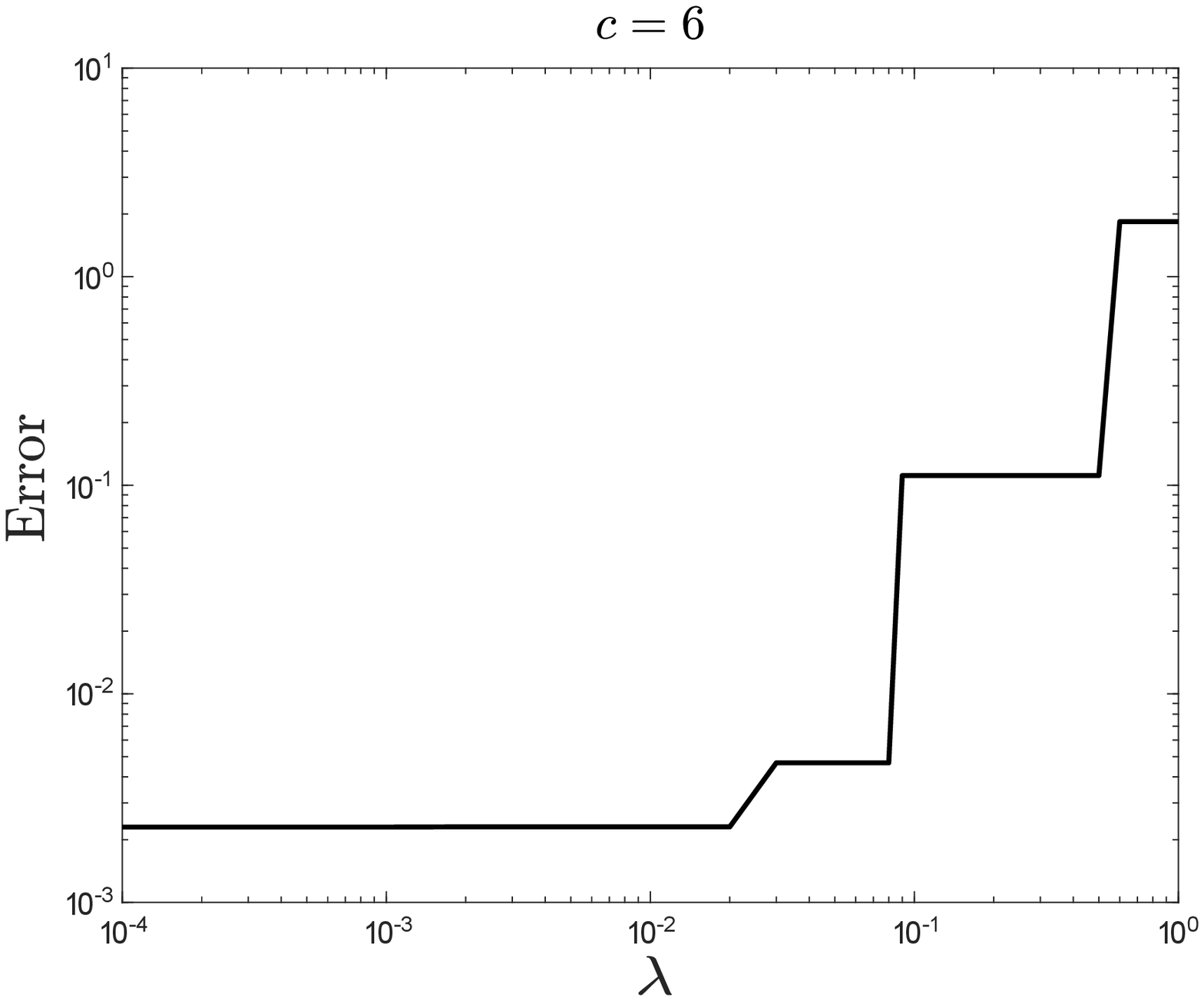} \quad
\includegraphics[width=0.35\textwidth,height = 0.32\textwidth]{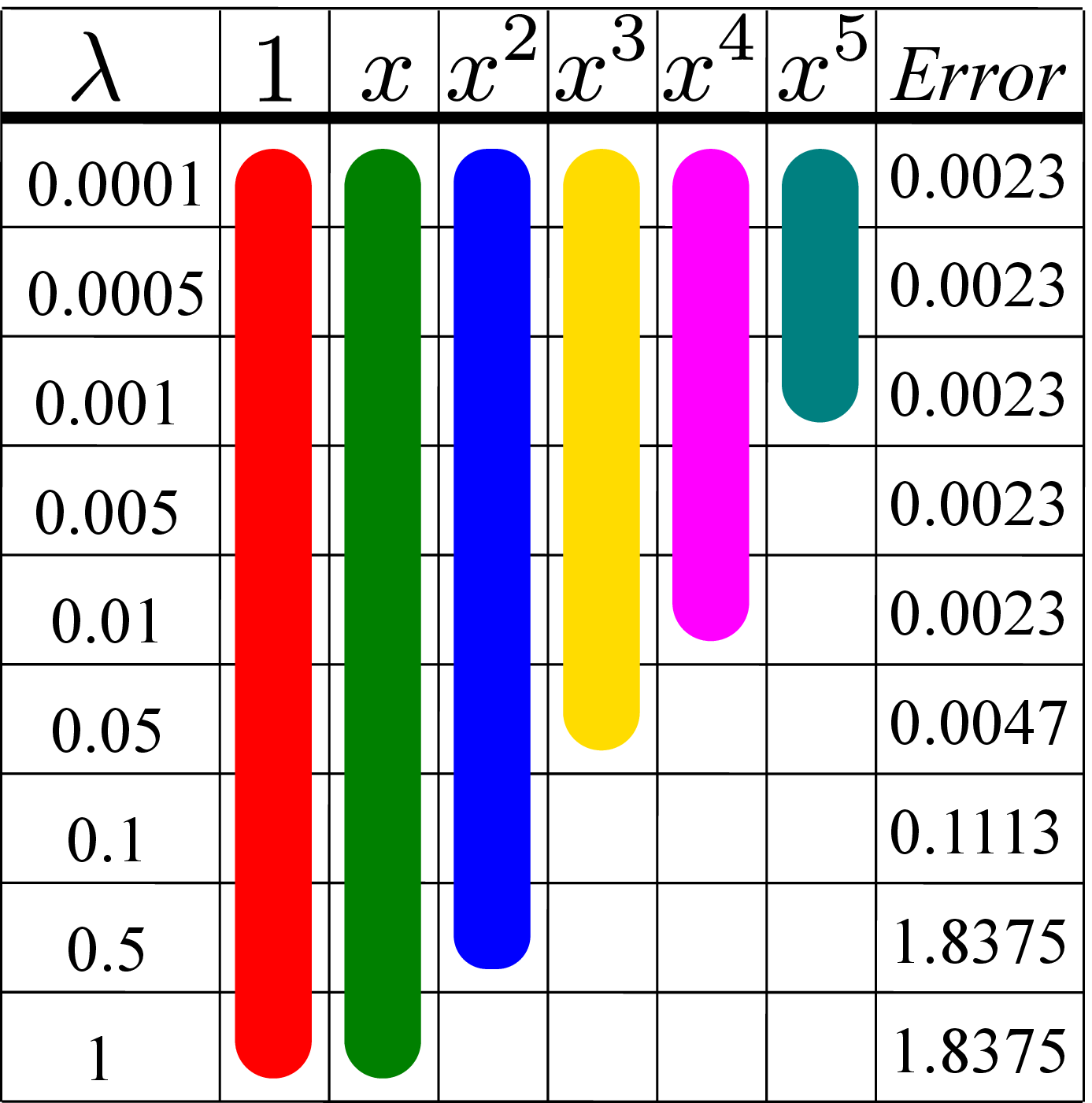}  \\
\includegraphics[width=0.45\textwidth]{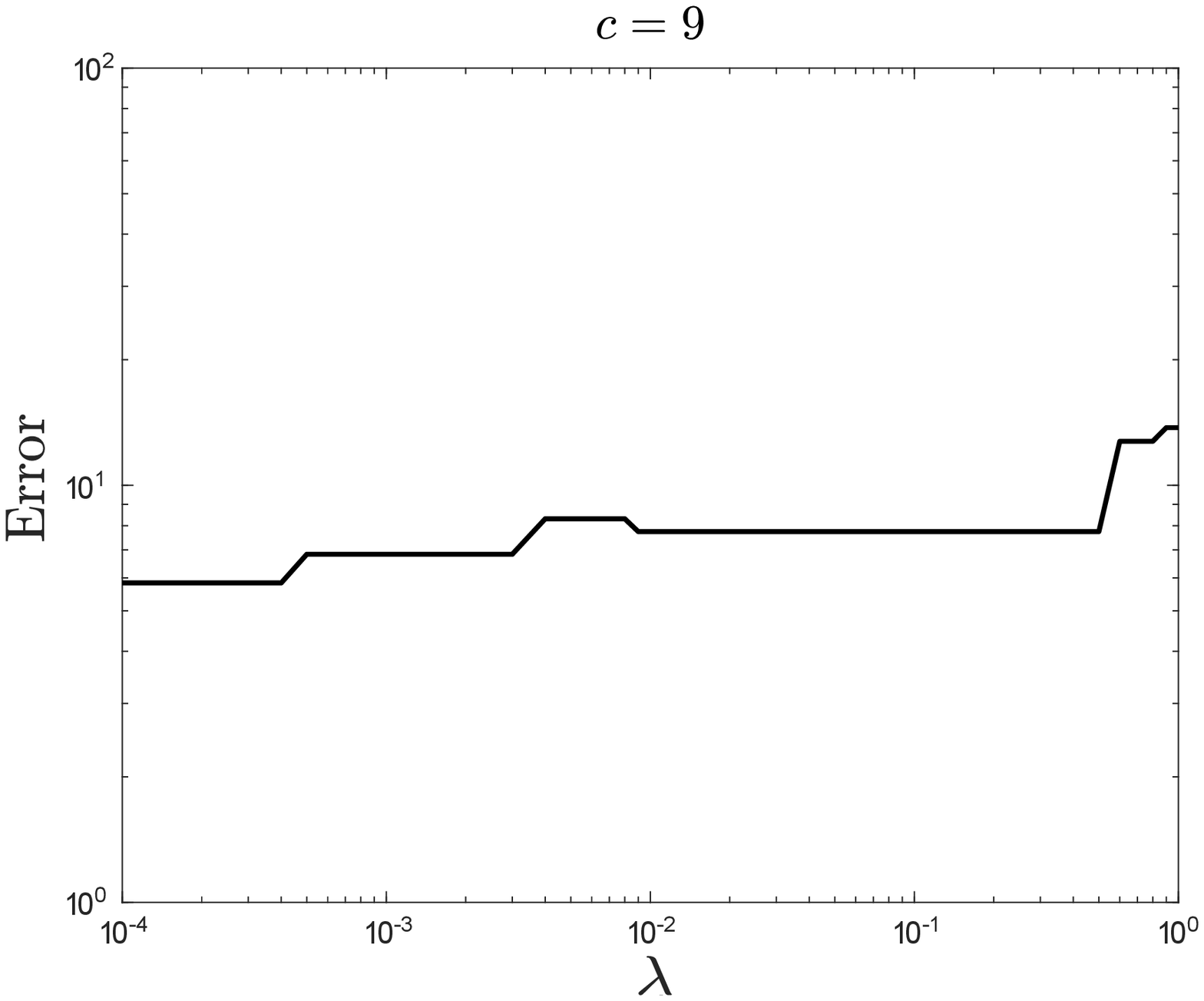} \quad
\includegraphics[width=0.35\textwidth,height = 0.32\textwidth]{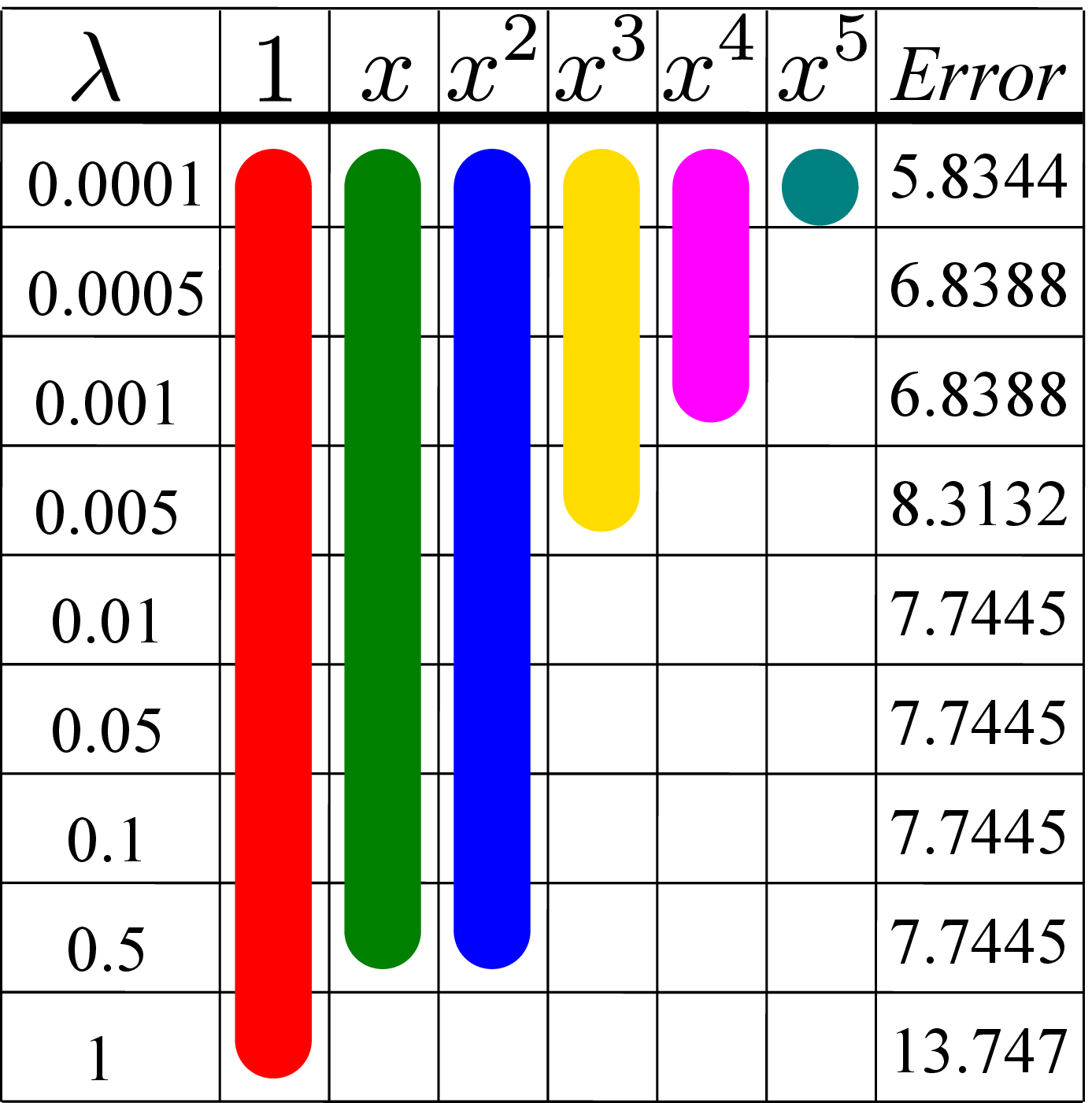} \\
\includegraphics[width=0.45\textwidth]{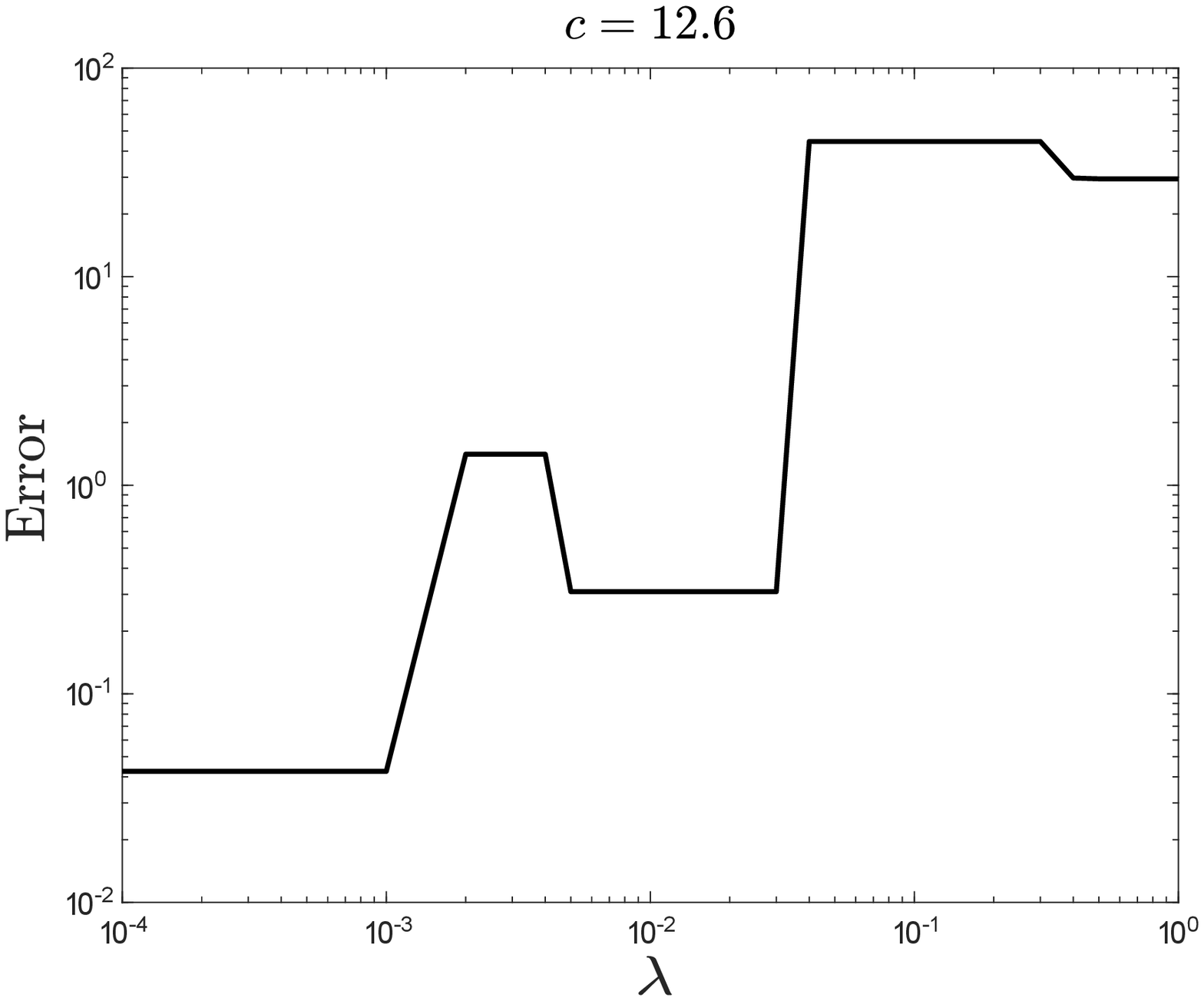} \quad
\includegraphics[width=0.35\textwidth,height = 0.32\textwidth]{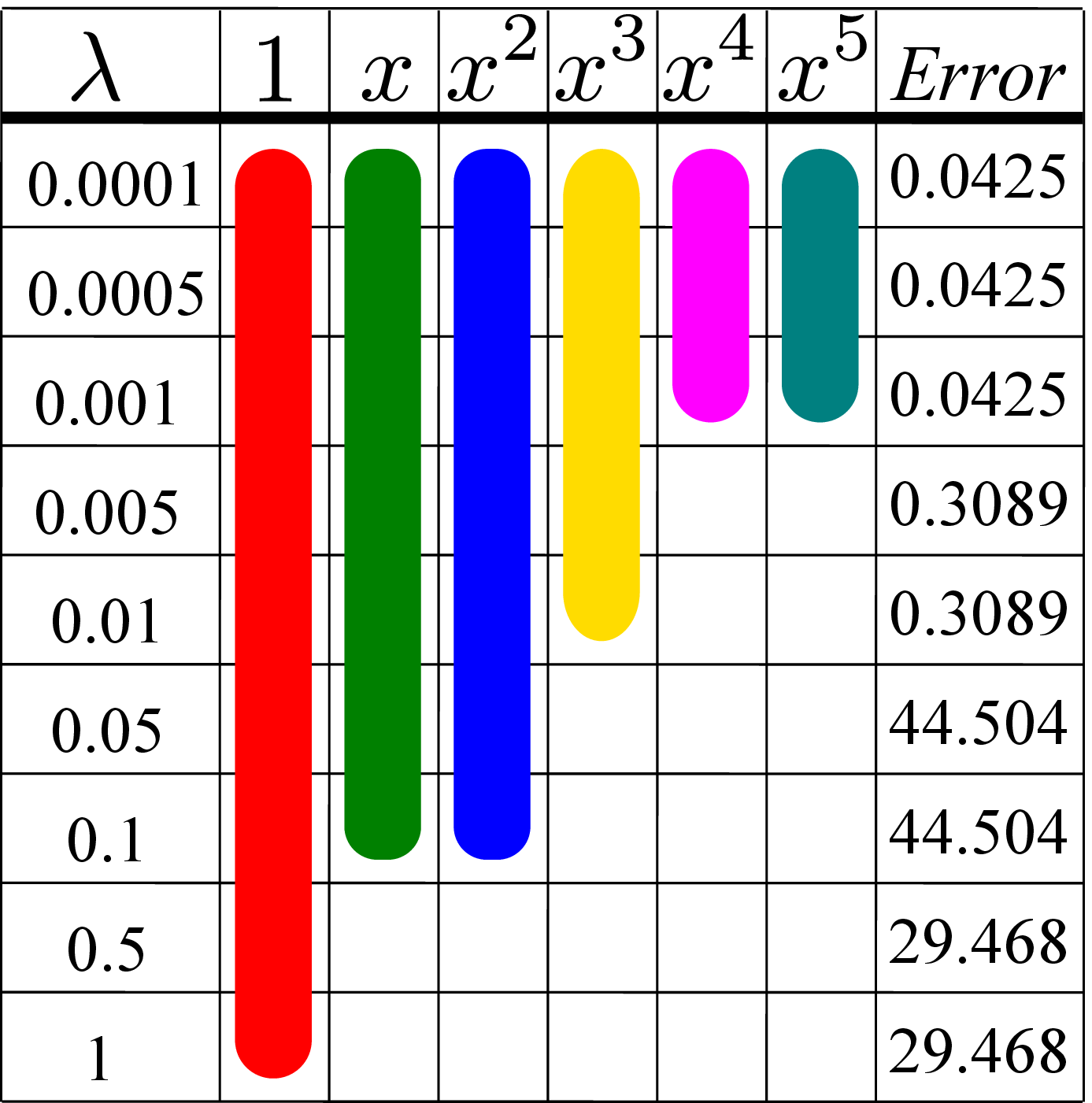}  \\
\includegraphics[width=0.45\textwidth]{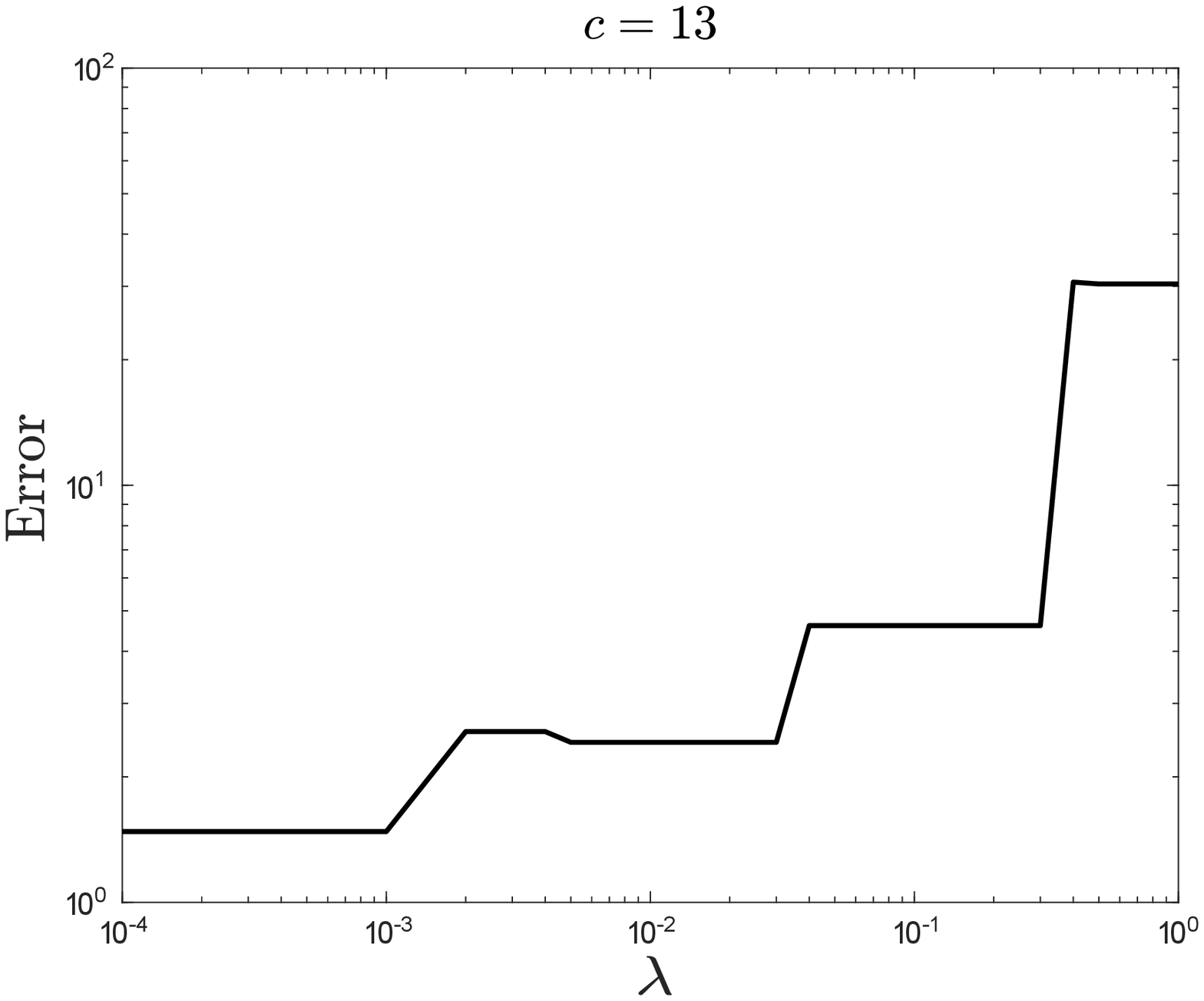} \quad
\includegraphics[width=0.35\textwidth,height = 0.32\textwidth]{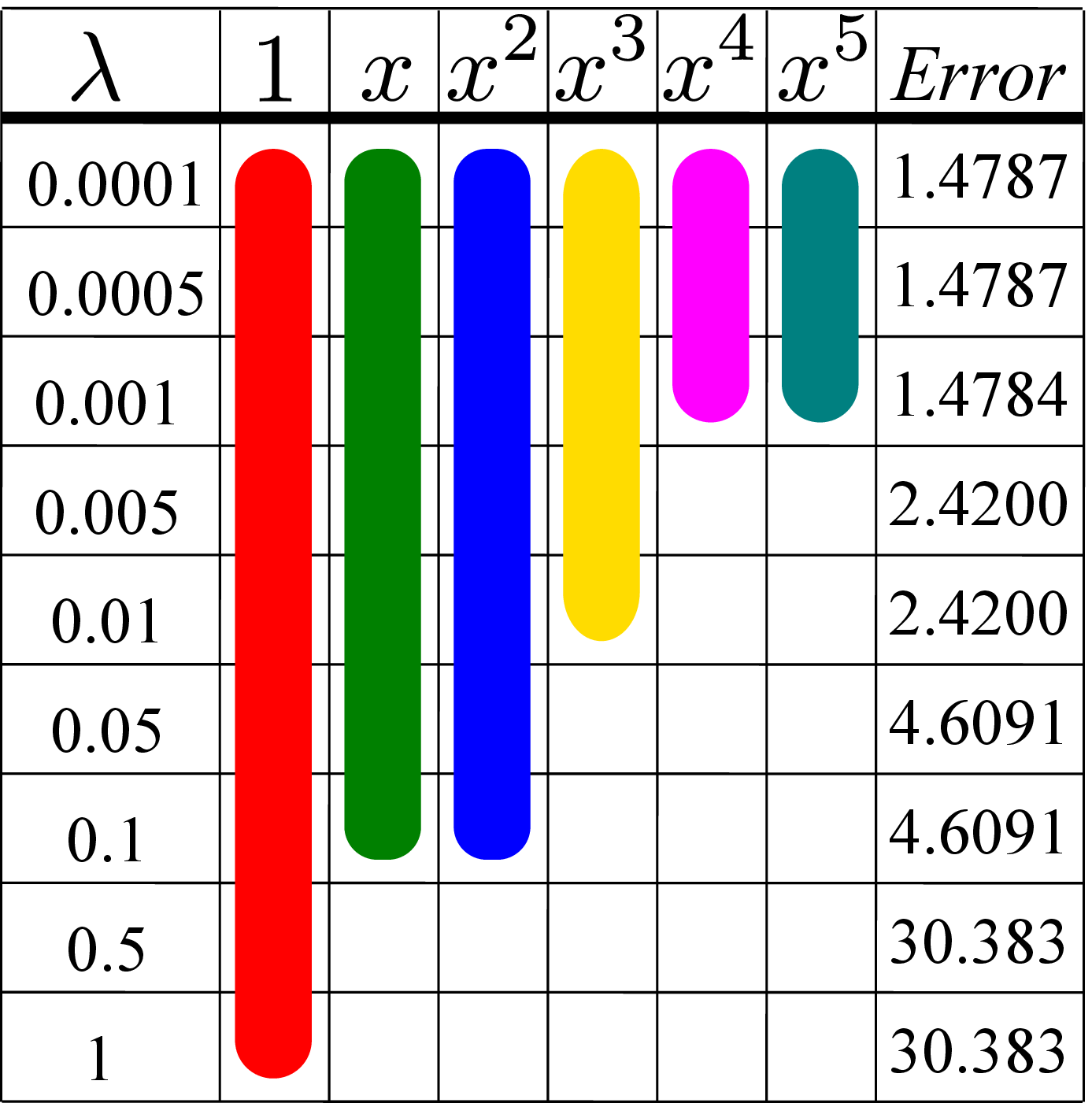} 
\caption{Error and sparsity of the obtained Poincar\'e mapping as a function of the sparsity parameter $\lambda$. (left) A log-log plot of the $\ell^2$ error between the training data and the simulated data with the same initial condition. (right) A visual depiction of the sparsity of the right-hand-side as a function of $\lambda$, along with the resulting error.}
\label{fig:RosslerError}
\end{figure}  

In Figure~\ref{fig:RosslerError} we provide different visualizations of the sensitivity of the method with respect to the sparsity parameter, $\lambda$. On the left we give log-log plots of the sparsity parameter against the $\ell^2$ distance between each iteration of the training data and the data generated by our mapping. On the right we provide a visualization of the values of the sparsity parameter which must be chosen to add each element of the library to the mapping, along with the error at each parameter value. As expected, decreasing the sparsity parameter overall decreases the error as well. Moreover, one can see that in the chaotic parameter regimes the error remains relatively large. We point out that this is again due to the fact that the mapping does not perfectly reproduce the chaotic dynamics but only obtains a seemingly chaotic attractor with the same banding structure.

\subsection{Spiral Waves} 

We now describe how our method can be applied to forecast solution behavior of partial differential equations. To illustrate this aspect of our method we consider a $\Lambda-\Omega$ reaction-diffusion system \cite{Murray}
\begin{equation}\label{LambdaOmega}
	\begin{split}
		u_t &= D\Delta u + \Lambda(A)u - \Omega(A)v, \\
		v_t &= D\Delta u + \Omega(A)u - \Lambda(A)v, \\
	\end{split}
\end{equation}
where $A = u^2 + v^2$, $\Lambda(A) = 1 - A^2$, $\Omega(A) = -\beta A^2$, and $D,\beta \geq 0$ are system parameters. Throughout this subsection we will fix $D = 0.1$ and $\beta = 1$ for illustration.  

\begin{figure} 
\center
\includegraphics[width=0.95\textwidth]{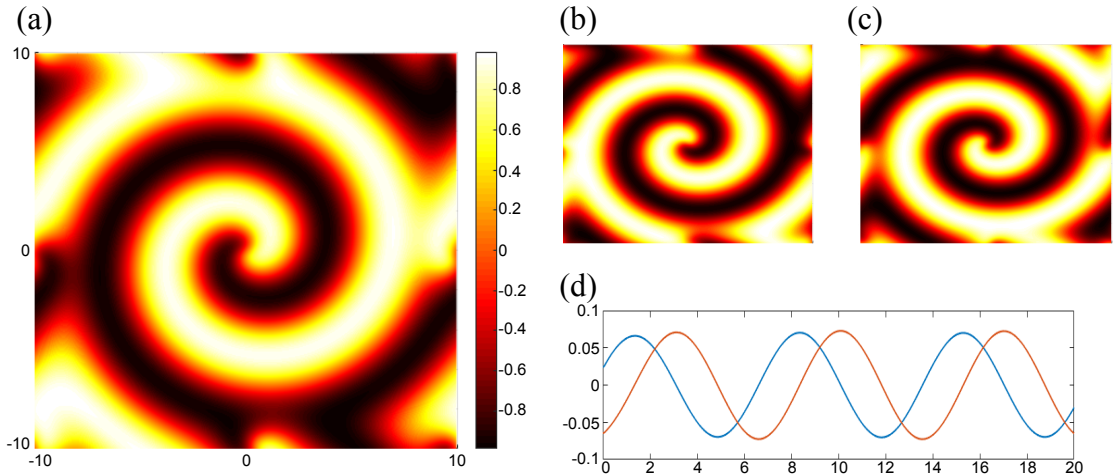} 
\caption{A spiral wave solution of (\ref{LambdaOmega}) with $D = 0.1$ and $\beta = 1$ along with its spatio-temporal singular value decomposition. (a) The $u$ component of the spiral wave solution of (\ref{LambdaOmega}) on the domain $x \in [-10,10]$ and $y \in [-10,10]$ with periodic boundary conditions. There are two dominant singular values associated to $u$ with left-singular vectors (b) and (c) associated to them. The temporal evolution of these left-singular vectors are nearly sinusoidal, with time series given in (d). The blue curve is associated to left-singular vector (b) whereas the orange curve is associated to left-singular vector (c).}
\label{fig:Spiral}
\end{figure}  

Reaction-diffusion systems of the type (\ref{LambdaOmega}) are well-known to exhibit rigidly rotating spiral waves \cite{Cohen,Ermentrout,Greenberg,Howard}. In Figure~\ref{fig:Spiral} we present a numerically obtained spiral wave solution of (\ref{LambdaOmega}) on the domain $x \in [-10,10]$ and $y \in [-10,10]$ with periodic boundary conditions. It was demonstrated in \cite{KoopmanPDE} that performing a space-time singular value decomposition on the numerically obtained solution data results in two dominant singular values for both $u$ and $v$ whose associated left-singular vectors simply oscillate sinusoidally in the temporal variable. This temporal oscillation of the dominant left-singular vectors gives the effect that the spiral wave is rotating about a single fixed point in space. These associated left-singular vectors are shown in Figure~\ref{fig:Spiral} along with the oscillatory right-singular vector associated to their temporal evolution. 

Due to the computational complexity of simulating accurate solutions to partial differential equations, we now propose that our work in this manuscript can be extended to provide accurate forecasting of nonlinear dynamics at a significantly decreased computational cost. In the context of the present example, we can imagine asking for accurate representations of the spiral wave long into the future without numerically integrating (\ref{LambdaOmega}) for a long time. To do this, we simulate the PDE on a short time interval with a small temporal step size and then reduce to the low-rank approximation coming from the spatio-temporal singular value decomposition. We may then extract snapshots of the temporal evolution from the right-singular vectors at regular time intervals to generate training data for our method. Upon applying our method with this training data, we will obtain a mapping which can be used to forecast snapshots of the spiral wave dynamics long into the future at regular time intervals. Our method would therefore take data about the current snapshot and provide a prediction of the following snapshot.    

\begin{figure} 
\center
\includegraphics[width=0.49\textwidth]{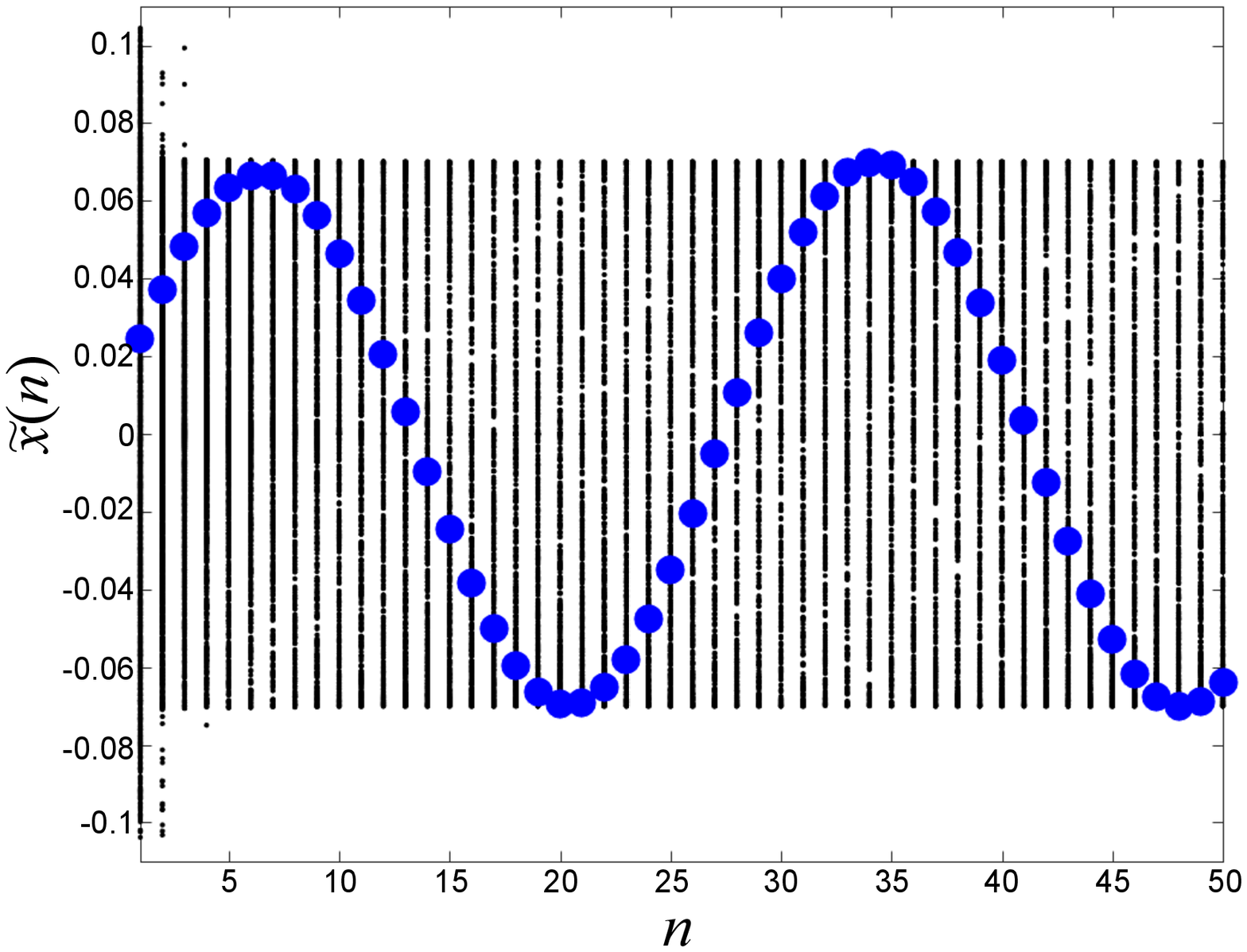} 
\includegraphics[width=0.49\textwidth]{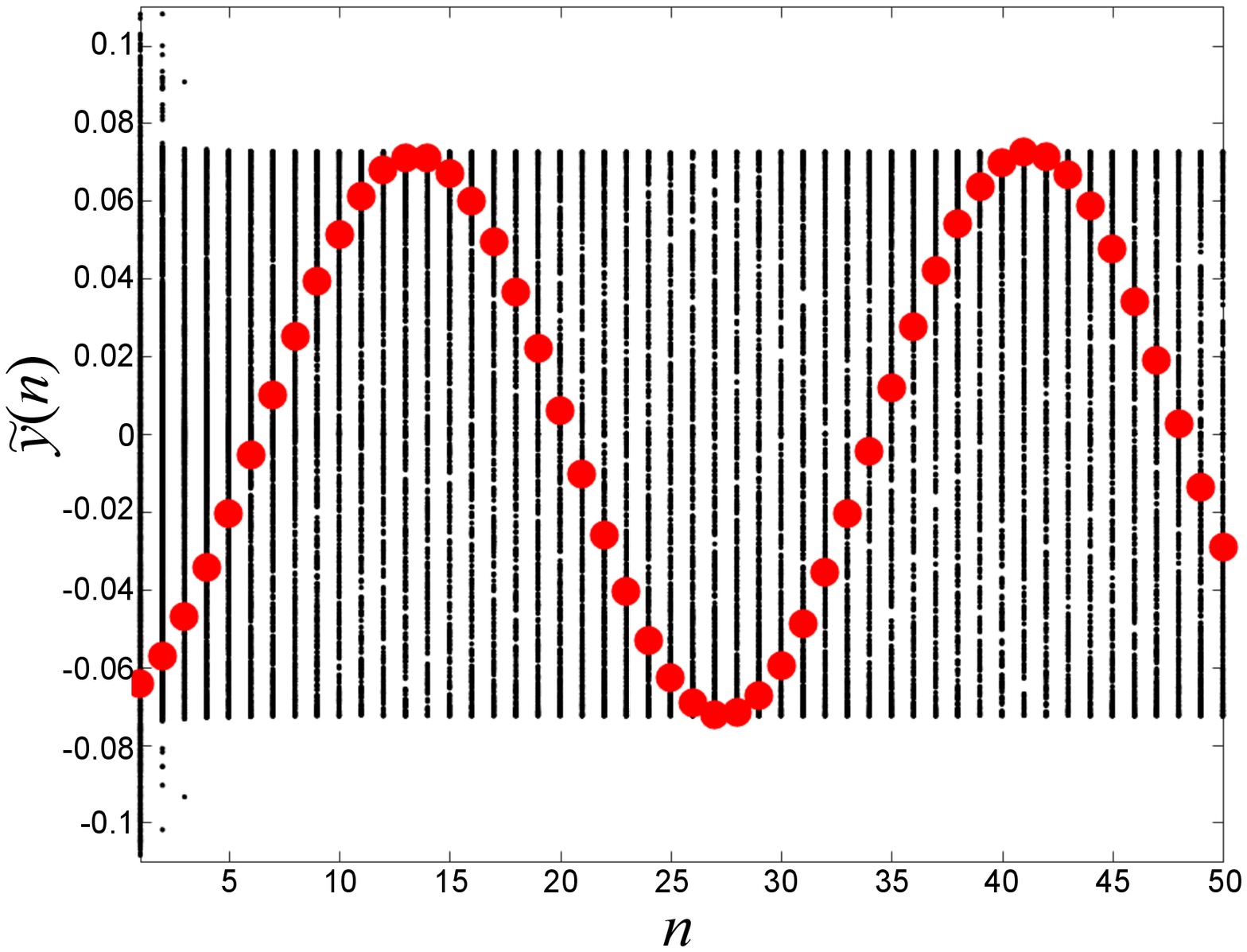}  
\caption{The first 50 iterations of the SINDy map resulting from sampling the time series Figure~\ref{fig:Spiral}(d) with $t \in \frac{1}{4}\mathbb{Z}\cap[0,20]$. Solutions converge quickly to a quasi-periodic attractor corresponding to the oscillation of the two dominant left-singular vectors in the spatio-temporal singular value decomposition of the spiral wave solution to (\ref{LambdaOmega}). Plotted on each of the images are also 1860 other trajectories (black) with random initial conditions.}
\label{fig:Spiral2}
\end{figure}  

To demonstrate this method, we simulated (\ref{LambdaOmega}) on the temporal domain $t \in [0,20]$ with a step-size $dt = 0.05$. Our snapshots are taken for each $t \in \frac{1}{4}\mathbb{Z}$ which corresponds to every fifth time-step. We note that the spiral wave dynamics are periodic in $t$, but unless the period of oscillation is commensurate with an element of $\frac{1}{4}\mathbb{Z}$ we would expect the training data to appear quasi-periodic. This is exactly what we see in Figure~\ref{fig:Spiral2}, where the large colored dots of blue and red represent the resulting SINDy map iterations with initial conditions taken to be the initial value of the temporal evolution of the left-singular vectors presented in Figure~\ref{fig:Spiral}(d). Here we again use the notation $(\tilde{x}(n),\tilde{y}(n))$ for the iterations of the SINDy map corresponding to the blue and orange curves of Figure~\ref{fig:Spiral2}, respectively. It follows that we may use this mapping to generate computationally cheap, accurate approximations of the spiral wave dynamics well beyond our limited timeframe $t \in [0,20]$ by simply iterating the map for as long as desired and using these generated values to form linear combinations of the dominant left-singular vectors. This in turn almost completely recovers the numerical spiral wave solution for arbitrarily large $t \geq 0$.    

We draw the readers attention to the similarity of this method to diffusion maps~\cite{coifman1,coifman2}, but note that our method presents both simulations on a dimensionally reduced manifold and a way of mapping these solutions on the manifold back (approximately) to the full solutions of the differential equation. In this way we are able to perform computationally cheap simulations of long-time dynamics of numerical solutions to partial differential equations. This is especially useful for multiscale systems where the fast scale can be parametrized by the Poincar\'e map and its discovered SINDy model.  Furthermore, this method allows one to analyze the effect of small deviations from initial conditions with ease using the obtained SINDy mapping, i.e. Monte-Carlo simulations can easily be done from a distribution of initial conditions. To exemplify this, in Figure~\ref{fig:Spiral2} we also present the evolution of 1860 other bounded trajectories of our obtained SINDy map with random initial conditions. This dimensional reduction therefore leads to a simple method of uncertainty quantification (UQ) that can be performed in an efficient manner that does not require simulating the full partial differential equation.  Indeed, one can compute the probability distribution on the Poincar\'e section from a given initial distribution to produce a proxy UQ metric at a substantial computational savings.

\section{Discussion}\label{sec:Discussion} 

In this manuscript we have seen the application of the SINDy method to discover Poincar\'e maps from data. Our methods were tested against a range of examples that included both ordinary and partial differential equations. As demonstrated in Section~\ref{sec:Applications}, these methods can be applied to provide long-time forecasting of the dynamics on and near invariant manifolds of differential equations. These methods have highlighted how the Poincar\'e map can be used to disambiguate and discover multiscale temporal dynamics, specifically the coarse-grained dynamics which results from fast-scale nonlinear dynamics.  

Potentially the most important aspect of applying this method is choosing the set of functions $\Thetav$ appropriately. In this work we simply took $\Thetav$ to contain polynomial functions of the state variables, but this is most likely not optimal. The reason for this is that even simple polynomial mappings can exhibit very complicated behavior on a set of measure zero while trajectories that remain in the complement of this set diverge. This can be observed in the simple quadratic mapping 
\[
	x_{n+1} = x_n^2 - \mu, \quad \mu > 2,
\] 
since almost every initial condition diverges but there is a Cantor set for which chaotic dynamics can be found. The simplicity of the quadratic mapping provides a warning that even simple one-dimensional nonlinear mappings can generate extremely complicated behavior for which forecasting of the related differential equation cannot be achieved.  

A potential avenue for future work could come from building in specific features of the approximate mapping $\tilde{\Piv}$, as was recently done for constructing continuous time conservation laws from data \cite{Kaiser}. For example, Poincar\'e mappings coming from dynamical systems with a conserved quantity would be expected to be measure preserving. Since minute perturbations can destroy a mappings measure-preserving property, it would not be expected that our method should provide a measure-preserving map due to potential numerical error introduced at each step in the method. Hence, augmenting the method herein to build in constraints would make the method more robust moving forward.   Additionally, for multiscale physics problems with a multitude of well separated timescales, one can imagine using the Poincar\'e map recursively to extract proxy models that parametrize each scale in a recursive and principled manner.

\end{document}